\documentclass[a4paper,12pt]{amsart}

\numberwithin{equation}{section}
\theoremstyle{plain}
\newtheorem{theorem}{\sc \bf Theorem}[section]
\newtheorem{lemma}[theorem]{\sc \bf Lemma}
\newtheorem{corollary}[theorem]{\sc \bf Corollary}
\newtheorem{proposition}[theorem]{\sc \bf Proposition}
\newtheorem{claim}{\sc \bf Claim}
\theoremstyle{definition}
\newtheorem{definition}[theorem]{\sc \bf Definition}
\newtheorem{remark}[theorem]{\sc \bf Remark}
\newtheorem{example}[theorem]{\sc \bf Example}

\usepackage{amssymb}
\usepackage{amsfonts}
\usepackage{amsmath}
\usepackage{enumerate}

\usepackage{tikz}

\newcommand{\K}{{\mathbb{K}}}
\newcommand{\R}{{\mathbb{R}}}
\newcommand{\C}{{\mathbb{C}}}

\newcommand{\ali}[1]{\begin{align*}#1\end{align*}}
\newcommand{\alil}[1]{\begin{align}#1\end{align}}

\newcommand{\si}{\sigma}
\newcommand{\Ph}{\varPhi}
\newcommand{\ph}{\varphi}
\newcommand{\phe}{\varphi(\epsilon,\cdot)}
\newcommand{\lam}{\lambda}
\newcommand{\e}{\epsilon}
\newcommand{\p}{{\prime}}

\newcommand{\BR}{\mathcal{B}}

\newcommand{\DR}{\mathcal{D}}
\newcommand{\ER}{\mathcal{E}}

\newcommand{\IR}{\mathcal{I}}

\newcommand{\LR}{\mathcal{L}}
\let\MR=\relax
\newcommand{\MR}{\mathcal{M}}

\newcommand{\OR}{\mathcal{O}}
\newcommand{\PR}{\mathcal{P}}
\newcommand{\QR}{\mathcal{Q}}
\newcommand{\RR}{\mathcal{R}}
\newcommand{\SR}{\mathcal{S}}
\newcommand{\TR}{\mathcal{T}}

\newcommand{\XR}{\mathcal{X}}
\newcommand{\YR}{\mathcal{Y}}

\newcommand{\tLR}{\tilde{\mathcal{L}}}

\newcommand{\ite}[1]{\begin{enumerate}[(1)]#1\end{enumerate}}
\let\prop=\relax

\newcommand{\prop}[1]{\begin{proposition}#1\end{proposition}}
\newcommand{\props}{\begin{proposition}}
\newcommand{\prope}{\end{proposition}}

\newcommand{\cors}{\begin{corollary}}
\newcommand{\core}{\end{corollary}}
\newcommand{\thm}[1]{\begin{theorem}#1\end{theorem}}
\newcommand{\thms}{\begin{theorem}}
\newcommand{\thme}{\end{theorem}}
\newcommand{\lem}[1]{\begin{lemma}#1\end{lemma}}
\newcommand{\lems}{\begin{lemma}}
\newcommand{\leme}{\end{lemma}}
\newcommand{\defi}[1]{\begin{definition}#1\end{definition}}
\newcommand{\defis}{\begin{definition}}
\newcommand{\defie}{\end{definition}}

\newcommand{\exams}{\begin{example}}
\newcommand{\exame}{\end{example}}
\newcommand{\rem}[1]{\begin{remark}\normalfont #1\end{remark}}
\newcommand{\cla}[1]{\begin{claim}#1\end{claim}}
\newcommand{\ncla}[1]{\setcounter{claim}{0}\begin{claim}#1\end{claim}}

\newcommand{\pros}{\begin{proof}}
\newcommand{\proe}{\end{proof}}
\newcommand{\prossq}{\begin{proof}}
\newcommand{\case}[1]{\begin{cases}#1\end{cases}}
\newcommand{\cd}{\cdot}
\newcommand{\up}{\upsilon}
\newcommand{\om}{\omega}
\newcommand{\ti}[1]{\tilde{#1}}

\newcommand{\di}{\displaystyle}

\newcommand{\qqqqqquad}{\qquad\qquad\qquad\qquad}
\newcommand{\qqqqqqquad}{\qquad\qquad\qquad\qquad\qquad}

\newcounter{constants}
\makeatletter
\def\addconst{
\addtocounter{constants}{1}
\def\@currentlabel{\arabic{constants}}
\@currentlabel
}
\makeatother

\newcommand{\adl}[1]{\addconst\label{c:#1}}
\newcommand{\adr}[1]{\ref{c:#1}}

\newcommand{\dist}{{\mathrm{dist}}}
\newcommand{\var}{\mathrm{var}}

\begin{document}
\title[Asymptotic perturbation]{General asymptotic perturbation theory in transfer operators}
\author[H. Tanaka]{Haruyoshi Tanaka}
\address{
{\rm Haruyoshi Tanaka}\\
Department of Mathematics and Statistics\\
Wakayama Medical University\\
580, Mikazura, Wakayama-city, Wakayama, 641-0011, Japan
}
\email{htanaka@wakayama-med.ac.jp}
\begin{abstract}
We study higher-order asymptotic expansions of eigenvalues in perturbed transfer operators, of the corresponding eigenfunctions and of the corresponding eigenvectors of the dual operators. In our main result, we give explicit expressions of these coefficients and these remainders under mild conditions of linear operators and of (there is even no norm) linear spaces. Moreover we investigate sufficient conditions for convergence of the remainders under weak conditions in the sense that the eigenvalue of perturbed operators does not have a uniform spectrum gap. As our main application, we apply our result to the asymptotic behaviour of Gibbs measures associated with the Hausdorff dimension of the limit set of  the so-called graph directed Markov systems. Moreover, we demonstrate the asymptotic behaviours of Ruelle transfer operators without a uniform spectral gap. In another example, we mention that the conditions of our main results are satisfied under the conditions of abstract asymptotic theory given by Gou\"ezel and Liverani and under mild conditions.
\end{abstract}
\maketitle
\tableofcontents
\section{Introduction}
Fix a nonnegative integer $n$. Put $\K=\R$ or $\C$. Let $\BR^{0},\BR^{1},\dots, \BR^{n+1}$ be linear spaces over $\K$ with $\BR^{0}\supset \BR^{1}\supset\cdots \supset \BR^{n+1}$ as linear subspaces. 
We consider an $n$-order asymptotic expansion of a perturbed linear operator $\LR(\e,\cd)\,:\,\BR^{0}\to \BR^{0}$ with a small parameter $\e\in (0,1)$:
\ali
{
\LR(\e,\cd)=\LR_{0}+\LR_{1}\e+\cdots+\LR_{n}\e^{n}+\tLR_{n}(\e,\cd)\e^{n} \text{ on }\BR^{n+1},
}
where the linear operator $\tLR_{n}(\e,\cd)$ is referred to as the remainder part of this expansion. Here we assume that each $\LR_{k}$ is a linear operator from $\BR^{i}$ to $\BR^{i-j}$ for each $0\leq j\leq n$ and $i=j,j+1,\dots, n$
(see the condition (I) in Section \ref{sec:rep_remainer}). Such a sequence of linear spaces and a sequence of linear operators are originally introduced by Gou\"ezel and Liverani \cite{GL}. We also consider the case that the operators $\LR$ and $\LR(\e,\cd)$ possess triplets $(\lam,h,\nu)$ and $(\lam(\e),h(\e,\cd),\nu(\e,\cd))$ in $\K\times \BR^{n+1}\times (\BR^{0})^{*}$ so that
\ali
{
&\LR h=\lam h,\ \ \LR^{*} \nu=\lam \nu,\ \ \nu(h)=1\\
&\LR(\e,\cd) h(\e,\cd)=\lam(\e) h(\e,\cd),\ \ \LR(\e,\cd)^{*} \nu(\e,\cd)=\lam(\e) \nu(\e,\cd),\ \ \nu(\e,h)\neq 0,\ \ \nu(h(\e,\cd))\neq 0
}
and $\lam$ is not an eigenvalue of the operator $\RR:=\LR-\lam (h\otimes \nu)$ acting on $\BR^{i}$ for $1\leq i\leq n$ (see the conditions (II), (III), (IV) in Section \ref{sec:rep_remainer} for detail).
\smallskip
\par
One of our main results is as follows: we give explicit expressions of each coefficients $\lam_{k}, \kappa_{k}, g_{k}$ and of the remainders $\ti{\lam}_{n}(\e)$, $\ti{\kappa}_{n}(\e,\cd)$, $\ti{g}_{n}(\e,\cd)$ of the expansions
\ali
{
\lam(\e)=&\lam+\lam_{1}\e+\cdots+\lam_{n}\e^{n}+\ti{\lam}_{n}(\e)\e^{n}\\
\frac{\nu(\e,f)}{\nu(\e,h)}=&\nu(f)+\kappa_{1}(f)\e+\cdots+\kappa_{n}(f)\e^{n}+\ti{\kappa}_{n}(\e,f)\e^{n}\\
\frac{h(\e,\cd)}{\nu(h(\e,\cd))}=&h+g_{1}\e+\cdots+g_{n}\e^{n}+\ti{g}_{n}(\e,\cd)\e^{n}
}
(Theorem \ref{th:rep_asymp_eval_nB_gene} and Theorem \ref{th:rep_asymp_eval_nB_gene2}). Note that the existence of $h(\e,\cd)$ is unnecessary to obtain the expansion of $\nu(\e,\cd)/\nu(\e,h)$ and the converse is similar. In order to show the asymptotic behaviours of these expansions (in a word, the remainder parts $\ti{\lam}_{n}(\e)$, $\ti{\kappa}_{n}(\e,f)$, $\ti{g}_{n}(\e,\cd)$ vanish as $\e\to 0$), we only have to give suitable norms and suitable conditions for vanishing the remainders of $\LR(\e,\cd)$. Consequently, we obtain sufficient conditions for $\ti{\lam}_{n}(\e)\to 0$ and $\ti{\kappa}_{n}(\e,f)\to 0$ for each $f$ (Theorem \ref{th:asymp_eval_nB_gene2}), for $\ti{\kappa}_{n}(\e,\cd)\to 0$ in a norm (Theorem \ref{th:asymp_eval_nB_gene3}) and for $\ti{g}_{n}(\e,\cd)\to 0$ (Theorem \ref{th:asymp_efunc_nB_gene} and Theorem \ref{th:asymp_efunc_nospecgap}). We note that in Theorem \ref{th:asymp_efunc_nospecgap}, we treat a sequence of Banach spaces depending on $\e$. We also stress that in Theorem \ref{th:asymp_eval_nB_gene2} and Theorem \ref{th:asymp_efunc_nospecgap}, we do not demand that the perturbed operator $\LR(\e,\cd)$ has a uniform spectral gap for $\lam$. Here we say that a perturbed operator $\LR(\e,\cd)$ has a uniform spectral gap for $\lam$ if the spectrum of $\LR(\e,\cd)$ in a fixed small open ball at $\lam$ consists only of $\lam(\e)$ for any small $\e>0$.
\smallskip
\par
The main motivation for the present work is to extend the result of asymptotic behaviour of the Gibbs measure associated with the Hausdorff dimension of the limit set of graph iterated function systems in the finite graph case \cite{T2011,T2016} to in the infinite graph case (Theorem \ref{th:asymp_Gibbs_GDMS}). We note that the Gibbs measure is displayed by using the corresponding eigenfunction and the corresponding eigenvector of the maximal positive eigenvalue of a suitable Ruelle transfer operator. Therefore our convergence results are useful for convergence of Gibbs measures. A second motivation is to remove the condition for the so-called uniform Lasota-Yorke type inequality in the perturbation theory of Ruelle transfer operators (see (GL.3) in Section \ref{sec:GL} for such a inequality). This inequality is important in the perturbation theory as the condition to provide with a uniform spectral gap. On the other hand, this inequality is not always necessary to obtain asymptotic behaviours of the topological pressure and the Gibbs measure for a perturbed potential. We will demonstrate asymptotically of the Gibbs measure for a suitable perturbed potential without a uniform spectral gap (therefore without a uniform Lasota-Yorke type inequality) (Theorem \ref{th:asymp_mue2} and Theorem \ref{th:ex_conv_Gibbs_nonspec}).
\smallskip
\par
To obtain asymptotic behaviours of eigenvalues and the corresponding eigenprojections, the asymptotic perturbation theory in \cite{Kato} is the most famous tool. The stability of spectrum of transfer operators (i.e. the case $n=0$) is also treated in \cite{KL, KL2}. Moreover, \cite{GL} gives an abstract asymptotic perturbation theory using a sequence of Banach spaces and applies this result to a perturbed Anosov diffeomorphism on smooth manifold. In these, a uniform spectral gap property of perturbed operators is imposed. This property yields the stability of the eigenprojection displayed by using the integration of the resolvent and consequently the isolated simple eigenvalue and the corresponding eigenvector continuously (or asymptotically) change. On the other hand, the spectral gap property is not always necessary to obtain asymptotic behaviours of the eigenvalue $\lam(\e)$ and of the corresponding eigenvector $\nu(\e,\cd)$. In fact, we gave in \cite{T2011} the asymptotic behaviour of the Perron eigenvalue and of the Perron eigenvector of a perturbed Ruelle operator for subshift of finite type without a spectral gap property. However, it was not found in the previous work whether the asymptotical behavior of the corresponding eigenfunction $h(\e,\cd)$ is obtained without such a gap property.
The main results (Theorem \ref{th:rep_asymp_eval_nB_gene} and Theorem \ref{th:rep_asymp_eval_nB_gene2}) of the present work are proved by generalizing the technique for inductively giving coefficients of expansions of the Ruelle operator in \cite{T2011}. By investigating conditions for vanishing remainders $\ti{\lam}_{n}(\e)$, $\ti{\kappa}_{n}(\e,\cd)$ and $\ti{g}_{n}(\e,\cd)$, it turn out that the result of the asymptotical behaviour can be applied to a wide class of bounded linear operators (including Ruelle operator with countable state space). In particular, we obtain a sufficient conditions for these vanishing without a spectral gap property (Theorem \ref{th:asymp_eval_nB_gene2} and Theorem \ref{th:asymp_efunc_nospecgap}). Note that our method has the characteristic of avoiding a perturbation of resolvent, but it is a simpler technique. The formulation for main results is inspired by \cite{GL}. 
\smallskip
\par
In application, we extend the previous results \cite{T2011,T2016} of the asymptotic behaviours of the thermodynamical features (the topological pressure and the Gibbs measure) for a perturbed potential defined on subshift of finite type to on countable Markov shift (Theorem \ref{th:asymp_pphe}, Theorem \ref{th:asymp_pphe_he}). Furthermore, we apply the behaviour of the Gibbs measure associated with the dimension of the limit set of a perturbed graph directed Markov system to these results (Theorem \ref{th:asymp_Gibbs_GDMS}). This result shall be important in a future work for the multifractal analysis of such a system using the asymptotic analysis. On the other hand, in the study of the perturbed Gibbs measure, Theorem \ref{th:asymp_mue2} enables the perturbed potential $\phe$ such that the spectral gap of the Ruelle operator $\phe$ vanishes (Section \ref{sec:ex_withoutSG} for example). In another application, we mention that the conditions of our main results are satisfied under the conditions of abstract asymptotic theory given by \cite{GL} and under mild conditions (Theorem \ref{th:GLtoGAPT}). As a result, it is found to be able to apply the present study to the perturbation problems of various dynamical systems (e.g. \cite{Baladi2,BC}). 
\smallskip
\par
In Section \ref{sec:rep_remainer}, we state and prove our main results for asymptotic expansions of the eigenvalue $\lam(\e)$, of the corresponding eigenfunction $h(\e,\cd)$, and of the corresponding eigenvector $\nu(\e,\cd)$ of the dual operators.
Additionally, we give in Section \ref{sec:conv_remainder} sufficient conditions for vanishing the remainder parts of $\lam(\e)$, of $\nu(\e,f)/\nu(\e,h)$ and of $h(\e,\cd)/\nu(h(\e,\cd))$. In particular, we introduce a notion of weak boundedness in Definition \ref{def:weakbdd} which is a generalization of a properly of the unperturbed resolvent given in \cite{KL, GL} and which is used in convergence of eigenfunction (Theorem \ref{th:asymp_efunc_nB_gene}). We will give in Proposition \ref{prop:prop_weakbdd} useful necessary and sufficient conditions for weak boundedness. On one hand, Theorem \ref{th:asymp_efunc_nospecgap} treats the case that each Banach space depends on the parameter $\e>0$. This enables us to deal with perturbations of eigenfunctions without a uniform spectral gap. In Section \ref{sec:coefn=012}, we examine the coefficients and the remainders under the cases $n=0,1,2$ of expansions given in Section \ref{sec:rep_remainer}. As application in Section \ref{ex:CMS}, we describe asymptotic expansions and give sufficient conditions for vanishing remainder parts for thermodynamic features of perturbed potentials defined on countable Markov shift. In Section \ref{sec:GDMS_Gibbs}, we apply asymptotic perturbations of graph directed Markov systems to the results of Section \ref{ex:CMS}. The section \ref{sec:ex_withoutSG} presents behaviours of thermodynamic features without a uniform spectral gap of Ruelle operators with subshift of finite type. We mention in Section \ref{sec:GL} a relation between the abstract perturbation theory of \cite{GL} and our results. In appendix \ref{sec:thermo}, we recall the notions of countable Markov shifts, of thermodynamic formalism, and of Ruelle operators which are needed in Section \ref{ex:CMS}. 
\medskip
\\
{\it Acknowledgment.}\ 
This study was supported by JSPS KAKENHI Grant Number 20K03636.
\section{Main results}
In the former subsection, we formulate asymptotic perturbations of linear operators acting on linear spaces and give exact expressions of coefficients and remainders of expansions of eigenvalues, eigenfunctions, and eigenvectors of the dual operator. In the later subsection, by using the expressions of remainders, we will introduce sufficient conditions for convergence of the remainders which are useful to prove our applications.
\subsection{Representation of the remainder parts of asymptotic perturbations under abstract setting}\label{sec:rep_remainer}
Fix a nonnegative integer $n$. Put $\K=\R$ or $\C$. Let $\BR^{0},\BR^{1},\dots, \BR^{n+1}$ be linear spaces over $\K$ with $\BR^{0}\supset \BR^{1}\supset\cdots \supset \BR^{n+1}$ as linear subspaces. At first we introduce linear operators $\LR=\LR_{0},\LR_{1},\dots, \LR_{n}, \LR(\e,\cd)$ satisfying the following (I)-(III):
\ite
{
\item[(I)] Each $\LR_{j}$ is a linear operator from $\BR^{i}$ to $\BR^{i-j}$ for each $j=0,1,\dots, n$ and $i=j,j+1,\dots, n$.
\item[(II)] The operator $\LR\,:\,\BR^{0}\to \BR^{0}$ has the decomposition $\LR=\lam (h\otimes \nu)+\RR$ such that (i) $(\lam,h,\nu)\in \K\times \BR^{n+1}\times (\BR^{0})^{*}$ satisfies $\lam\neq 0$, $\LR^{*}\nu =\lam \nu$, $\LR h=\lam h$ and $\nu(h)=1$, (ii) there exists a linear subspace $\DR^{0}$ with $\BR^{1}\subset \DR^{0}\subset \BR^{0}$ such that the inverse $(\RR-\lam\IR)^{-1}$ from $\DR^{0}$ to $\BR^{0}$ is well-defined and $(\RR-\lam\IR)^{-1}\BR^{i}\subset \BR^{i}$ for each $1\leq i\leq n+1$.
\item[(III)] The parametrized linear operator $\LR(\e,\cd)\,:\,\BR^{0}\to \BR^{0}$ with parameter $\e\in (0,1)$ has a pair $(\lam(\e),\nu(\e,\cd))\in \K\times (\BR^{0})^{*}$ satisfying that $\LR(\e,\cd)^{*}\nu(\e,\cd)=\lam(\e)\nu(\e,\cd)$ and $\nu(\e,h)\neq 0$.
}
Here the operator $h\otimes \nu\,:\,\BR^{0}\to \BR^{n+1}$ is given by $h\otimes \nu(f)=\nu(f)h$. 
\rem
{
If $\lam$ is not eigenvalue of $\RR\,:\,\BR^{0}\to \BR^{0}$, then we may put $\DR^{0}=(\RR-\lam\IR)\BR^{0}$. In general, $\BR^{1}\varsubsetneq\DR^{0}\varsubsetneq \BR^{0}$ is satisfied (see \cite[Proposition 3.3]{T2011} for example).
}
For later convenience, we put
\alil
{
&\tLR_{j}(\e,\cd):=(\LR(\e,\cd)-\LR-\LR_{1}\e-\cdots-\LR_{j}\e^{j})/\e^{j}\ \ (j=0,1,\dots, n)\\
&\PR:=h\otimes \nu,\ \RR_{\lam}:=(\RR-\lam\IR)^{-1},\ \QR(\e,\cd):=h\otimes\kappa(\e,\cd)-\IR,\  \kappa(\e,\cd):=\nu(\e,\cd)/\nu(\e,h)\\
&\DR^{k}:=\{f\in \DR^{0}\,:\,\RR_{\lam}f\in \BR^{k+1},\ (\lam_{l}\IR-\LR_{l})\RR_{\lam}f\in \DR^{n-l}\quad (l=1,\dots, k)\}.\label{eq:Dk=...}
}
Note that $\DR^{k}$ is defined by induction and includes $\BR^{k+1}$. 
We give explicit formulations of asymptotic expansions of $\lam(\e)$, $g(\e,\cd)$ and $\kappa(\e,\cd)$: 
\thms
\label{th:rep_asymp_eval_nB_gene}
Assume that the conditions (I)-(III) are satisfied. Then we have the asymptotic expansions
\alil
{
\lam(\e)=&\lam+\lam_{1}\e+\cdots+\lam_{n}\e^{n}+\ti{\lam}_{n}(\e)\e^{n}\label{eq:le=...}\\
\kappa(\e,f)=&\nu(f)+\kappa_{1}(f)\e+\cdots+\kappa_{n}(f)\e^{n}+\ti{\kappa}_{n}(\e,f)\e^{n}\quad \text{for }f\in \DR^{n}\label{eq:kaef=...}
}
by putting
\alil
{
\lam_{k}=&\nu(\LR_{k}h)+\sum_{j=2}^{k}\sum_{i_{1},\dots,i_{j}\geq 1\,:\,\atop{i_{1}+\cdots+i_{j}=k}}\nu((\lam_{i_{1}}\IR-\LR_{i_{1}})\RR_{\lam}\cdots (\lam_{i_{j-1}}\IR-\LR_{i_{j-1}})\RR_{\lam}\LR_{i_{j}}h)\label{eq:lamk=...}\\
\kappa_{k}(f)=&\sum_{j=1}^{k}\sum_{i_{1},\dots,i_{j}\geq 1\,:\,\atop{i_{1}+\cdots+i_{j}=k}}\nu((\lam_{i_{1}}\IR-\LR_{i_{1}})\RR_{\lam}\cdots (\lam_{i_{j}}\IR-\LR_{i_{j}})\RR_{\lam}f)\label{eq:kapakf=}
}
for each $1\leq k\leq n$. Moreover, the remainders $\ti{\lam}_{k}(\e)$, $\ti{\kappa}_{k}(\e,\cd)$ and $\ti{g}_{k}(\e,\cd)$ have the forms
\alil
{
\ti{\lam}_{k}(\e)=&\kappa(\e,\tLR_{k}(\e,h))+\sum_{i_{1}\geq 0, i_{2}\geq 1\,:\,i_{1}+i_{2}=k}\kappa(\e,\tLR_{i_{1}}(\e,\cd)\QR(\e,\cd)\RR_{\lam}\LR_{i_{2}}h)\label{eq:tlamke=}\\
&+\sum_{p=3}^{k+1}\sum_{i_{1}\geq 0,i_{2},\dots,i_{p}\geq 1\atop{i_{1}+\dots+i_{p}=k}}\kappa(\e,\tLR_{i_{1}}(\e,\cd)\QR(\e,\cd)\RR_{\lam}(\prod_{j=1}^{p-1}(\LR_{i_{j}}\QR(\e,\cd)\RR_{\lam}+\lam_{i_{j}}\RR_{\lam}))\LR_{i_{p}}h)\nonumber\\
\ti{\kappa}_{k}(\e,f)=&\kappa(\e,\tLR_{k}(\e,\cd)\QR(\e,\cd)\RR_{\lam} f)\label{eq:tkakef=}\\
+&\sum_{p=2}^{k+1}\sum_{i_{1}\geq 0,i_{2},\dots,i_{p}\geq 1\,:\,\atop{i_{1}+\cdots+i_{p}=k}}\kappa(\e,\tLR_{i_{1}}(\e,\cd)\QR(\e,\cd)\RR_{\lam}\prod_{j=2}^{p}(\LR_{i_{j}}\QR(\e,\cd)\RR_{\lam}+\lam_{i_{j}}\RR_{\lam})f)\nonumber
}
for each $0\leq k\leq n$ and $\e>0$. 
\rem
{
Assume that the conditions (I)-(III) are satisfied. The following inductive definitions of coefficients and remainders are useful in our applications:
\alil
{
\lam_{k}=&\sum_{j=1}^{k}\kappa_{k-j}(\LR_{j}h),\qquad
\kappa_{k}(f)=\sum_{j=1}^{k}\kappa_{k-j}((\lam_{j}\IR-\LR_{j})\RR_{\lam} f)\label{eq:lamk=kappak=...in}
}
for each $k=1,2,\dots, n$ inductively, where $\kappa_{0}=\nu$ and $g_{0}=h$, and
\alil
{
\ti{\lam}_{k}(\e)=&\kappa(\e,\tLR_{k}(\e,h))+\sum_{l=1}^{k}\ti{\kappa}_{k-l}(\e,\LR_{l}h)\label{eq:tlamke=...in}\\
\ti{\kappa}_{k}(\e,f)=&\kappa(\e,\tLR_{k}(\e,\cd)\QR(\e,\cd)\RR_{\lam} f)+\sum_{l=1}^{k}\ti{\kappa}_{k-l}(\e,(\LR_{l}\QR(\e,\cd)\RR_{\lam}+\lam_{l}\RR_{\lam}) f)\label{eq:tkapke=...in}
}
for each $0\leq k\leq n$. Note that these inductive definitions imply the equations (\ref{eq:lamk=...})-(\ref{eq:tkakef=}).
}
\pros
We define $\lam_{k}$ and $\kappa_{k}$ by (\ref{eq:lamk=kappak=...in}). These are well-defined by $h\in \BR^{n+1}\subset \BR^{k+1}$, $\LR_{j}h\in \BR^{k+1-j}\subset \DR^{k-j}$ for $1\leq j\leq k$, and by $(\lam_{j}\IR-\LR_{j})\RR_{\lam}\DR^{k}\subset (\lam_{j}\IR-\LR_{j})\BR^{k+1}\subset \BR^{k+1}\cap \BR^{k+1-j}=\BR^{k+1-j}\subset \DR^{k-j}$. Similarity, $\QR(\e,\cd)\RR_{\lam}$ is well-defined on $\DR^{0}$. We start with the following:
\ncla
{\label{cla:I-P=(L-l)R_l}
$\IR-\PR=(\LR-\lam\IR)\RR_{\lam}$ on $\DR^{0}$. 
}
Indeed, let $x\in \DR^{0}$ and take $y\in \RR_{\lam}\DR^{0}$ satisfying $x=(\RR-\lam\IR)y$. We have
\ali
{
(\IR-\PR)x=&(\IR-\PR)(\RR-\lam\IR)y=
(\LR-\lam\IR)y=(\LR-\lam\IR)\RR_{\lam}x.
}
\cla
{
The equations (\ref{eq:le=...}) and (\ref{eq:kaef=...}) hold for $k=0$.
}
Indeed, consider the equation
\ali
{
\nu(\e,(\lam(\e)-\lam)h)&=\nu(\e,(\LR(\e,\cd)-\LR)h)
}
using $\lam(\e)\nu(\e,h)=\nu(\e,\LR(\e,h))$ and $\LR h=\lam h$. Then we have $\lam(\e)-\lam=\kappa(\e,(\LR(\e,\cd)-\LR)h)=\kappa(\e,\tLR_{0}(\e,h))=\ti{\lam}_{0}(\e)$. Moreover, we obtain that
for $f\in \DR^{0}$
\ali
{
\kappa(\e,f)-\nu(f)=&\kappa(\e, f-h\nu(f))\qquad (\because \kappa(\e,h)\equiv 1)\\
=&\kappa(\e,(\IR-\PR)f)\\
=&\kappa(\e,(\LR-\lam \IR)\RR_{\lam} f)\qquad (\because \text{Claim }\ref{cla:I-P=(L-l)R_l})\\
=&\kappa(\e,(\ti{\lam}_{0}(\e)\IR-\tLR_{0}(\e,\cd))\RR_{\lam} f)\\
=&\kappa(\e,\cd)\tLR_{0}(\e,\cd)\QR(\e,\cd)\RR_{\lam} f=\ti{\kappa}_{0}(\e,f).
}
Thus we obtain the claim.
\cla
{\label{cla:asymp_lam_kappa}
If the equations (\ref{eq:le=...}) and (\ref{eq:kaef=...}) hold for $k^\p=0,1,\dots, k-1$ with $k\geq 1$, then so do for $k^\p=k$.
}
We have
\ali
{
\frac{\lam(\e)-\lam-\sum_{i=1}^{k}\lam_{i}\e^{i}}{\e^{k}}=&\kappa(\e,\frac{\lam(\e)-\lam}{\e^{k}}h)-\sum_{i=1}^{k}\frac{\lam_{i}\e^{i}}{\e^k}\\
=&\kappa(\e,\frac{\LR(\e,\cd)-\LR-\sum_{l=1}^{k}\LR_{l}\e^{l}}{\e^{k}}h)+\kappa(\e,\sum_{l=1}^{k}\frac{\LR_{l}\e^{l}}{\e^{k}}h)-\sum_{i=1}^{k}\frac{\lam_{i}\e^{i}}{\e^{k}}\\
=&\kappa(\e,\tLR_{k}(\e,h))+\sum_{l=1}^{k}\kappa(\e,\LR_{l}h)\e^{l-k}-\sum_{i=1}^{k}\sum_{j=1}^{i}\kappa_{i-j}(\LR_{j}h)\e^{i-k}\\
=&\kappa(\e,\tLR_{k}(\e,h))+\sum_{l=1}^{k}\kappa(\e,\LR_{l}h)\e^{l-k}-\sum_{j=1}^{k}\sum_{i=0}^{k-j}\kappa_{i}(\LR_{j}h)\e^{i+j-k}=\ti{\lam}_{k}(\e),
}
where the last equation uses (\ref{eq:tlamke=...in}), Moreover,
\ali
{
&\frac{\kappa(\e,f)-\nu(f)-\sum_{i=1}^{k}\kappa_{i}(f)\e^{i}}{\e^{k}}\\
=&\frac{\kappa(\e,(\IR-\PR)f)}{\e^{k}}-\sum_{i=1}^{k}\frac{\kappa_{i}(f)\e^{i}}{\e^{k}}\\
=&\frac{\kappa(\e,(\LR-\LR(\e,\cd)+\lam(\e)\IR-\lam\IR)\RR_{\lam} f)}{\e^{k}}-\sum_{i=1}^{k}\sum_{j=1}^{i}\kappa_{i-j}((\lam_{j}\IR-\LR_{j})\RR_{\lam} f)\e^{i-k}\\
=&\frac{\lam(\e)-\lam-\sum_{l=1}^{k}\lam_{l}\e^{l}}{\e^{k}}\kappa(\e,\RR_{\lam} f)-\kappa(\e,\frac{\LR(\e,\cd)-\LR-\sum_{l=1}^{k}\LR_{l}\e^{l}}{\e^{k}}\RR_{\lam} f)\\
&+\sum_{l=1}^{k}\kappa(\e,(\lam_{l}\IR-\LR_{l})\RR_{\lam} f) \e^{l-k}-\sum_{j=1}^{k}\sum_{i=j}^{k}\kappa_{i-j}((\lam_{j}\IR-\LR_{j})\RR_{\lam} f)\e^{i-k}\\
=&\kappa(\e,(\ti{\lam}_{k}(\e)-\tLR_{k}(\e,\cd))\RR_{\lam} f)+\sum_{l=1}^{k}\ti{\kappa}_{k-l}(\e,(\lam_{l}\IR-\LR_{l})\RR_{\lam} f)\\
=&\kappa(\e,\cd)\tLR_{k}(\e,\cd)\QR(\e,\cd)\RR_{\lam} f+\sum_{l=1}^{k}\ti{\kappa}_{k-l}(\e,\cd)(\LR_{l}\QR(\e,\cd)\RR_{\lam}+\lam_{l}\RR_{\lam}) f=\ti{\kappa}_{k}(\e,f).
}
Hence the proof is complete.
\proe
\thme
We next state the asymptotic expansions of eigenfunctions. We introduce the following condition:
\ite
{
\item[(IV)] The linear operator $\LR(\e,\cd)\,:\,\BR^{n+1}\to \BR^{0}$ has a pair $(\lam(\e),h(\e,\cd))\in \K\times \BR^{n+1}$ satisfying that $\LR(\e,h(\e,\cd))=\lam(\e)h(\e,\cd)$ and $\nu(h(\e,\cd))\neq 0$.
}
We set
\alil
{
g(\e,\cd):=h(\e,\cd)/\nu(h(\e,\cd)),\ \ \SR:=\RR_{\lam}(\IR-\PR),\ \ \TR(\e,\cd)=g(\e,\cd)\otimes \nu-\IR.\label{eq:ge=,S=,Te=}
}
We obtain the next result:
\thms
\label{th:rep_asymp_eval_nB_gene2}
Assume that the conditions (I),(II) and (IV) are satisfied. Then we have the asymptotic expansions
\alil
{
\lam(\e)=&\lam+\lam_{1}\e+\cdots+\lam_{n}\e^{n}+\ti{\lam}_{n}(\e)\e^{n}\label{eq:le=...2}\\
g(\e,\cd)=&h+g_{1}\e+\cdots+g_{n}\e^{n}+\ti{g}_{n}(\e)\e^{n}\label{eq:ge=...}
}
by putting
\alil
{
\lam_{k}=&\nu(\LR_{k}h)+
\sum_{j=2}^{k}\sum_{\underset{i_{1}+\cdots+i_{j}=k}{i_{1},\dots,i_{j}\geq 1\,:\,}}\nu(\LR_{i_{1}}\SR(\lam_{i_{2}}\IR-\LR_{i_{2}})\cdots \SR(\lam_{i_{j}}\IR-\LR_{i_{j}})h)\label{eq:lamk=_mth2}\\
g_{k}=&\sum_{j=1}^{k}\sum_{\underset{i_{1}+\cdots+i_{j}=k}{i_{1},\dots,i_{j}\geq 1\,:\,}}\SR(\lam_{i_{1}}\IR-\LR_{i_{1}})\cdots \SR(\lam_{i_{j}}\IR-\LR_{i_{j}})h\label{eq:gk=2}\\
\ti{\lam}_{k}(\e)=&\nu(\tLR_{k}(\e,g(\e,\cd)))+\sum_{i_{1}\geq 1,i_{2}\geq 0\,:\,i_{1}+i_{2}=k}\nu(\LR_{i_{1}}\SR\TR(\e,\cd)\tLR_{i_{2}}(\e,g(\e,\cd)))\label{eq:tlke=...2}\\
&+\sum_{p=3}^{k+1}\sum_{i_{1},\dots, i_{p-1}\geq 1,i_{p}\geq 0\,:\,\atop{i_{1}+\cdots+i_{p}=k}}\nu(\LR_{i_{1}}(\prod_{j=2}^{p-1}(\SR\TR(\e,\cd)\LR_{i_{j}}+\lam_{i_{j}}\SR))\SR\TR(\e,\cd)\tLR_{i_{p}}(\e,g(\e,\cd)))\nonumber\\
\ti{g}_{k}(\e,\cd)=&\SR\TR(\e,\cd)\tLR_{k}(\e,g(\e,\cd))\label{eq:tgke=...in2}\\
+&\sum_{p=2}^{k+1}\sum_{i_{1},\dots, i_{p-1}\geq 1,i_{p}\geq 0\,:\,\atop{i_{1}+\cdots+i_{p}=k}}(\prod_{j=1}^{p-1}(\SR\TR(\e,\cd)\LR_{i_{j}}+\lam_{i_{j}}\SR))\SR\TR(\e,\cd)\tLR_{i_{p}}(\e,g(\e,\cd))\nonumber
}
for each $0\leq k\leq n$ and $\e>0$.
\thme
\rem
{
The representations (\ref{eq:tlamke=}) and (\ref{eq:tlke=...2}) of $\ti{\lam}_{k}(\e)$ are different.
}
\rem
{
Assume that the conditions (I),(II) and (IV) are satisfied. We see
\alil
{
\lam_{k}=&\sum_{j=1}^{k}\nu(\LR_{j}g_{k-j}),\qquad
g_{k}=\sum_{j=1}^{k}\SR(\lam_{j}\IR-\LR_{j})g_{k-j}\label{eq:lamk=gk=...}
}
for each $k=1,2,\dots, n$ inductively, where $g_{0}=h$, and
\alil
{
\ti{\lam}_{k}(\e)=&\nu(\tLR_{k}(\e,g(\e,\cd)))+\sum_{l=1}^{k}\nu(\LR_{l}\ti{g}_{k-l}(\e,\cd))\label{eq:tlamke=...in2}\\
\ti{g}_{k}(\e,\cd)=&\SR\TR(\e,\cd)\tLR_{k}(\e,g(\e,\cd))+\sum_{l=1}^{k}(\SR\TR(\e,\cd)\LR_{l}+\SR \lam_{l})\ti{g}_{k-l}(\e,\cd)\label{eq:tgke=...in}
}
for each $0\leq k\leq n$. We remark that The inductive definitions yield (\ref{eq:lamk=_mth2})-(\ref{eq:tgke=...in2}).
}
\pros
We define $\lam_{k}$ and $g_{k}$ by (\ref{eq:lamk=gk=...}). These are well-defined by $h\in \BR^{n+1}\subset \BR^{k+1}$, $\LR_{j}h\in \BR^{k+1-j}\subset \DR^{k-j}$ for $1\leq j\leq k$, and by $(\lam_{j}\IR-\LR_{j})\RR_{\lam}\DR^{k}\subset (\lam_{j}\IR-\LR_{j})\BR^{k+1}\subset \BR^{k+1}\cap \BR^{k+1-j}=\BR^{k+1-j}\subset \DR^{k-j}$. We also remark that the operator $\SR$ is well-defined on $\DR^{0}$ and $\SR\BR^{i}\subset \BR^{i}$ for $1\leq i\leq n+1$. Similarity, $\QR(\e,\cd)\RR_{\lam}$ and $\SR \TR(\e,\cd)$ are well-defined on $\DR^{0}$. We start with the following:
\ncla
{\label{cla:I-P=(L-l)R_l2}
$\IR-\PR=\SR(\LR-\lam\IR)$ on $\RR_{\lam}\DR^{0}$. 
}
Choose any $x\in \RR_{\lam}\DR^{0}$. By $h\in \BR^{1}\subset \RR_{\lam}\DR^{0}$, we see $\PR x\in\RR_{\lam}\DR^{0}$ and then $y=(\RR-\IR\lam)(\IR-\PR)x$ is in $\DR^{0}$. We obtain
\ali
{
(\IR-\PR)(\LR-\lam\IR)x=&(\IR-\PR)(\RR-\lam\IR)x=y\in \DR^{0}.
}
Thus we see $\SR(\LR-\lam\IR)x=(\RR-\lam\IR)^{-1}y=(\IR-\PR)x$.
\cla
{
The equations (\ref{eq:le=...2}) and (\ref{eq:ge=...}) hold for $k=0$.
}
Indeed, the equation $\nu((\lam(\e)-\lam)g(\e,\cd))=\nu((\LR(\e,\cd)-\LR)g(\e,\cd))$ implies $\lam(\e)-\lam=\nu(\tLR_{0}(\e,g(\e,\cd)))=\ti{\lam}_{0}(\e)$. Furthermore,
\alil
{
g(\e,\cd)-h=&(\IR-\PR)g(\e,\cd)\qquad(\because \nu(g(\e,\cd))\equiv 1)\label{eq:ge-h=...}\\
=&\SR(\LR-\lam\IR) g(\e,\cd)\ \ (\because g(\e,\cd)\in \BR^{1}\subset \RR_{\lam}\DR^{0})\nonumber\\
=&\SR(-\tLR_{0}(\e,\cd)+\ti{\lam}_{0}(\e)\IR)g(\e,\cd)=\ti{g}_{0}(\e,\cd).\nonumber
}
Therefore we obtain the claim.
\cla
{
If The equations (\ref{eq:le=...2}) and (\ref{eq:ge=...}) hold for $k^\p=0,1,\dots, k-1$ with $k\geq 1$, then so do for $k^\p=k$.
}
The equation (\ref{eq:le=...2}) follows from a similar argument of Claim \ref{cla:asymp_lam_kappa} in Theorem \ref{th:rep_asymp_eval_nB_gene}.
We have
\ali
{
\frac{g(\e,\cd)-\sum_{i=0}^{k}g_{i}\e^{i}}{\e^{k}}=&\frac{(\IR-\PR)g(\e,\cd)}{\e^{k}}-\sum_{i=1}^{k}g_{i}\e^{i-k}\\
=&\SR(\frac{\lam(\e)-\lam}{\e^{k}}-\frac{\LR(\e,\cd)-\LR}{\e^{k}})g(\e,\cd)-\sum_{i=1}^{k}g_{i}\e^{i-k}\\
=&\SR((\ti{\lam}_{k}(\e)\IR-\tLR_{k}(\e,\cd))g(\e,\cd))\\
&+\SR(\sum_{l=1}^{k}\lam_{l}\e^{l-k}g(\e,\cd)-\sum_{l=1}^{k}\LR_{l}g(\e,\cd)\e^{l-k}-\sum_{j=1}^{k}\sum_{i=j}^{k}(\lam_{j}\IR-\LR_{j})g_{i-j}\e^{i-k})\\
=&\SR(\ti{\lam}_{k}(\e)\IR-\tLR_{k}(\e,\cd))g(\e,\cd)+\sum_{l=1}^{k}\SR(\lam_{l}\IR-\LR_{l})\ti{g}_{k-l}(\e,\cd)\\
=&\SR\TR(\e,\cd)\tLR_{k}(\e,g(\e,\cd))+\sum_{l=1}^{k}\SR\TR(\e,\cd)\LR_{l}\ti{g}_{k-l}(\e,\cd)+\sum_{l=1}^{k}\SR(\lam_{l}\ti{g}_{k-l}(\e,\cd))\\
=&\ti{g}_{k}(\e,\cd).
}
Thus we get (\ref{eq:ge=...}). Hence the proof is obtained.
\proe
Finally we give asymptotic expansions of $\nu(\e,\cd)$ and $h(\e,\cd)$. 
In addition to the conditions (I)-(IV), we impose the following:
\ite
{
\item[(V)] There exists an element $1_{\BR}\in \BR^{n+1}$ such that $\nu(1_{\BR})=\nu(\e,1_{\BR})=1$ for any $\e>0$.
\item[(VI)] $\nu(\e,h(\e,\cd))= 1$ for any $\e>0$.
}
We have the following:
\cors
\label{cor:asymp_nue_he}
Assume that the conditions (I)-(VI) are satisfied. Then 
\alil
{
\nu(\e,\cd)=&\nu+\nu_{1}\e+\cdots+\nu_{n}\e^{n}+\ti{\nu}_{n}(\e,\cd)\e^{n}\label{eq:nue=...}\\
h(\e,\cd)=&h+h_{1}\e+\cdots+h_{n}\e^{n}+\ti{h}_{n}(\e,\cd)\e^{n}\label{eq:he=...}
}
by putting
\alil
{
\nu_{k}=&\kappa_{k}+\sum_{i=1}^{k}\Big(\sum_{l=1}^{i}\sum_{j_{1},\dots, j_{l}\geq 1\,:\,\atop{j_{1}+\cdots+j_{l}=i}}(-1)^{l}\kappa_{j_{1}}(1_{\BR})\cdots \kappa_{j_{l}}(1_{\BR})\Big)\kappa_{k-i},\label{eq:nuk=...}\\
h_{k}=&g_{k}+\sum_{i=1}^{k}\Big(\sum_{l=1}^{i}\sum_{j_{1},\dots, j_{l}\geq 1\,:\,\atop{j_{1}+\cdots+j_{l}=i}}(-1)^{l}c(j_{1})\cdots c(j_{i})\Big)g_{k-i}\\
\ti{\nu}_{n}(\e,\cd)=&\sum_{k=1}^{n}\ti{a}_{n-k}(\e)\kappa_{k}+\nu(\e,h)\ti{\kappa}_{n}(\e,\cd)\\
\ti{h}_{n}(\e,\cd)=&\sum_{i=1}^{n}g_{i}\ti{c}_{n-i}(\e)+\ti{g}_{n}(\e,\cd)\nu(h(\e,\cd)),\label{eq:thne=...}
}
where $c(k)=\sum_{i=1}^{k}\nu_{i}(g_{k-i})$,
\alil
{
\ti{a}_{n}(\e)=&\sum_{k=2}^{n}\sum_{p=1}^{k}\sum_{j_{1},\dots, j_{p-1}\geq 1,j_{p}\geq 0\,:\,\atop{j_{1}+\cdots+j_{p}=n}}(-1)^{k}\kappa_{j_{1}}(1_{\BR})\cdots\kappa_{j_{p-1}}(1_{\BR})\ti{\kappa}_{j_{p}}(\e,1)\ti{\kappa}_{0}(\e,1_{\BR})^{k-p}\quad \text{and}\\
\ti{c}_{n}(\e)=&\sum_{k=1}^{n}\nu_{k}(\ti{g}_{n-k}(\e,\cd))+\ti{\nu}_{n}(\e,g(\e,\cd))\label{eq:tcne=...}\\
=&\sum_{k=0}^{n}\ti{\nu}_{n-k}(\e,g_{k})+\nu(\e,\ti{g}_{n}(\e,\cd)).\label{eq:tcne=...2}
}
In particular, $1/\nu(h(\e,\cd))=\sum_{k=0}^{n}c(k)\e^{k}+\ti{c}_{n}(\e)\e^{n}$ holds.
\core
\pros
When we put $f=1_{\BR}$ in the expansion (\ref{eq:kaef=...}) of $\kappa(\e,f)$, we have the expression
$1/\nu(\e,h)=1+\sum_{k=1}^{n}\kappa_{k}(1_{\BR})\e^{k}+\ti{\kappa}_{n}(\e,1_{\BR})\e^{n}$
and $\nu(1_{\BR})=1$. By using the Taylor expansion of the map $F\,:\,x\mapsto 1/x$ at $x=1$,
\ali
{
&\nu(\e,h)=\frac{1}{1/\nu(\e,h)}=1+\sum_{i=1}^{n}a_{i}\e^{i}+\ti{a}_{n}(\e)\e^{n},
}
where
$a_{i}=\sum_{k=1}^{i}\sum_{j_{1},\dots, j_{k}\geq 1\,:\,\atop{j_{1}+\cdots+j_{k}=i}}(-1)^{k}\kappa_{j_{1}}(1)\cdots \kappa_{j_{k}}(1)$.
Thus we have
\ali
{
\nu(\e,f)=&\kappa(\e,h)\nu(\e,h)\\
=&(\sum_{k=0}^{n}\kappa_{k}(f)\e^{k}+\ti{\kappa}_{n}(\e,f)\e^{n})\nu(\e,h)\\
=&\sum_{k=0}^{n}\sum_{i=0}^{n-k}a_{i}\kappa_{k}\e^{i+k}+\sum_{k=1}^{n}\kappa_{k}(f)\ti{a}_{n-k}(\e)\e^{n}+\nu(\e,h)\ti{\kappa}_{n}(\e,f)\e^{n}
}
and (\ref{eq:nue=...}) is fulfilled. On the other hand, we remark the expansion
\ali
{
\nu\Big(\e,\frac{h(\e,\cd)}{\nu(h(\e,\cd))}\Big)=&\frac{1}{\nu(h(\e,\cd))}\\
=&1+\sum_{k=1}^{n}\sum_{i,j\geq 0\,:\,i+j=k}\nu_{i}(g_{j})\e^{k}+\Big(\sum_{k=1}^{n}\nu_{k}(\ti{g}_{n-k}(\e,\cd))+\ti{\nu}_{n}(\e,\frac{h(\e,\cd)}{\nu(h(\e,\cd))})\Big)\e^{n}\\
=&1+\sum_{k=1}^{n}c_{k}\e^{k}+\ti{c}_{n}(\e)\e^{n}.
}
Thus
\ali
{
\nu(h(\e,\cd))=&1+\sum_{k=1}^{n}d_{k}\e^{k}+\ti{d}_{n}(\e)\e^{n}
}
with $d_{i}=\sum_{k=1}^{i}\sum_{j_{1},\dots, j_{k}\geq 1\,:\,\atop{j_{1}+\cdots+j_{k}=i}}(-1)^{k}c_{j_{1}}\cdots c_{j_{k}}$.
Consequently we obtain the expansion (\ref{eq:he=...}) of $h(\e,\cd)=\nu(h(\e,\cd))g(\e,\cd)$.
\proe
\subsection{Convergence of the remainder parts of asymptotic perturbations}\label{sec:conv_remainder}
For Banach spaces $\XR,\YR$, denoted by $\LR(\XR,\YR)$ the set of all bounded linear operators from $\XR$ to $\YR$ endowed with the operator norm $\|\LR\|_{\XR\to\YR}=\sup_{f\in \XR\,:\,\|f\|_{\XR}\leq 1}\|\LR f\|_{\YR}$ for $\LR\in \LR(\XR,\YR)$. If $\XR=\YR$ then we may write $\|\cd\|_{\XR\to\XR}=\|\cd\|_{\XR}$ for simple. We begin with the sufficient condition for convergence of the remainder of eigenvector $\nu(\e,\cd)$.
\thms
\label{th:asymp_eval_nB_gene2}
Assume that the conditions (I)-(III) are satisfied. We also assume that
\ite
{
\item[(a)] $\di\limsup_{\e\to 0}|\kappa(\e,f)|<+\infty$ for any $f\in \BR^{0}$ and;
\item[(b)] $\kappa(\e,\tLR_{j}(\e,f))\to 0$ as $\e\to 0$ for any $f\in \DR^{j}$ for $j=0,1,\dots, n$.
}
Then
$\ti{\lam}_{n}(\e)\to 0$ and $\ti{\kappa}_{n}(\e,f)\to 0$ for each $f\in \DR^{n}$. Moreover, if the condition (V) is also satisfied and $\sup_{\e>0}|\nu(\e,h)|<\infty$, then $\ti{\nu}_{n}(\e,f)\to 0$ for each $f\in \DR^{n}$.
\thme
\pros
Recall the forms (\ref{eq:tlamke=...in}) and (\ref{eq:tkapke=...in}). 
First we consider the case $k=0$. We have $\ti{\lam}_{0}(\e)=\kappa(\e,\tLR_{0}(\e,h))\to 0$ by the assumption (b). Moreover,
\ali
{
|\ti{\kappa}_{0}(\e,f)|=&|\kappa(\e,\tLR_{0}(\e,\cd)\QR(\e,\cd)\RR_{\lam}f)|\\
\leq &|\kappa(\e,\RR_{\lam}f)||\kappa(\e,\tLR_{0}(\e,\cd)h)|+ |\kappa(\e,\tLR_{0}(\e,\cd)f)|\to 0
}
for $f\in \DR^{0}$ by the conditions (a)(b) and by the fact $\RR_{\lam}f\in \BR^{1}$.
Next we check convergence in general case. To do this, we assume that $\ti{\lam}_{i}(\e)\to 0$ and $\ti{\kappa}_{i}(\e,f)\to 0$ for each $f\in \DR^{i}$ for $i=0,1,\dots, k-1$. In the form (\ref{eq:tlamke=...in}) of $\ti{\lam}_{k}(\e)$, we obtain $\ti{\lam}_{n}(\e)\to 0$ by the condition (b) with $\LR_{l}h\in \BR^{n+1-l}$. Furthermore, a similar argument yields convergence of $\kappa(\e,\tLR_{0}(\e,\cd)\QR(\e,\cd)\RR_{\lam}f)\to 0$. 
The assertion follows from $(\lam_{l}\IR-\LR_{l})\RR_{\lam}f\in \BR^{n-l}$ for $f\in \DR^{n}$ that
\ali
{
|\ti{\kappa}_{k-l}(\e,(\LR_{l}\QR(\e,\cd)\RR_{\lam}+\lam_{l}\RR_{\lam}) f)|\leq& |\kappa(\e,\RR_{\lam}f)||\ti{\kappa}_{k-l}(\e,(\LR_{l}h))|+|\ti{\kappa}_{k-l}(\e,\lam_{l}\RR_{\lam}f-\LR_{l}\RR_{\lam}f)|\\
\to&0
}
as $\e\to 0$.
The last assertion follows from the form of the remainder $\ti{\nu}(\e,\cd)$ of $\nu(\e,\cd)$. Hence the proof is complete.
\proe
\rem
{
Theorem \ref{th:asymp_eval_nB_gene2} is a generalization of \cite[Theorem 2.1]{T2011}.
}
The following is used for showing the asymptotic expansions of Gibbs measures (see the proof of Theorem \ref{th:asymp_pphe_he}).
\thms
\label{th:asymp_eval_nB_gene3}
Assume that the conditions (I)-(III) are satisfied. 
We assume also the following (a)-(c):
\ite
{
\item[(a)] Each $\BR^{k}$ is a Banach space.
\item[(b)] $\nu\in \LR(\BR^{0},\K)$, $\LR_{j}\in \LR(\BR^{i},\BR^{i-j})$ for $j=0,1,\dots, n$ and $i=j,j+1,\dots, n+1$,  $\sup_{\e>0}\|\kappa(\e,\cd)\|_{\BR^{0}\to\K}<+\infty$ and $\RR_{\lam}\in \LR(\BR^{1},\BR^{0})\cap \LR(\BR^{j},\BR^{j})$ for $j=1,2,\dots, n+1$
\item[(c)] $\lim_{\e\to 0}\|\tLR_{j}(\e,\cd)\|_{\BR^{j+1}\to \BR^{0}}=0$ for each $j=0,1,\dots, n$.
}
Put $\|f\|_{k}:=\max_{1\leq i\leq k}\|f\|_{\BR^{i}}$ and $\BR^{k}_{1}:=(\BR^{k},\|\cd\|_{k})$. Then each $\kappa_{k}$ is in $\LR(\BR^{k}_{1},\K)$ and $\|\ti{\kappa}_{k}(\e,\cd)\|_{\BR_{1}^{k+1}\to\K}$ vanishes as $\e\to 0$ for each $0\leq k\leq n$. Moreover, if the condition (V) is satisfied and $\sup_{\e>0}|\nu(\e,h)|<+\infty$, then $\nu_{k}\in \LR(\BR^{k}_{1},\K)$ and $\|\ti{\nu}_{k}(\e,\cd)\|_{\BR_{1}^{k+1}\to\K}\to 0$.
\thme
\pros
We begin with the first assertion. Recall the form (\ref{eq:kapakf=}) of $\kappa_{k}$. Note that for $1\leq i\leq n$ and $0\leq k\leq n$, we have
\ali
{
\|(\lam_{i}\IR-\LR_{i})\RR_{\lam}f\|_{\BR^{k}}\leq& |\lam_{i}|\|\RR_{\lam}\|_{\BR^{k+\alpha}\to \BR^{k}}\|f\|_{\BR^{k+\alpha}}+\|\LR_{i}\|_{\BR^{i+l}\to \BR^{k}}\|\RR_{\lam}\|_{\BR^{k+i}\to \BR^{k+i}}\|f\|_{\BR^{k+i}},
}
where $\alpha=1$ if $k=0$ and $\alpha=0$ if $k>0$.
Therefore, for $i_{1},\dots, i_{j}\geq 1$, we have
\alil
{
\|(\prod_{l=1}^{j}(\lam_{i_{l}}\IR-\LR_{i_{l}})\RR_{\lam})f\|_{\BR^{k}}\leq&c \|f\|_{i_{1}+\cdots+i_{j}+k}\label{eq:bprodI-LRfBk}
}
for some constant $c>0$. Thus $\kappa_{k}$ is in $\LR(\BR^{k}_{1},\K)$.
\\
Finally we prove $\|\ti{\kappa}_{k}(\e,\cd)\|_{\BR_{1}^{k+1}\to\K}\to 0$ as $\e\to 0$. For $1\leq k\leq n$, $1\leq j\leq k$, $i_{1},\dots, i_{j}\geq 1$ with $i_{1}+\cdots+i_{j}=l$ and $f\in \BR^{k+1}$, we see
\ali
{
&|\kappa(\e,(\ti{\lam}_{k}(\e)\IR-\tLR_{k}(\e,\cd))\RR_{\lam}f)|\\
\leq& \|\kappa(\e,\cd)\|_{\BR^{0}\to\K}\{|\ti{\lam}_{k}(\e)|\|\RR_{\lam}\|_{\BR^{1}\to \BR^{0}}\|f\|_{\BR^{1}}
+\|\tLR_{k}(\e,\cd)\|_{\BR^{k+1}\to \BR^{0}}\|\RR_{\lam}\|_{\BR^{k+1}\to \BR^{k+1}}\|f\|_{\BR^{k+1}}\} 
}
and therefore
\ali
{
&|\kappa(\e,(\ti{\lam}_{i_{0}}(\e)\IR-\tLR_{i_{0}}(\e,\cd))\RR_{\lam}({\textstyle{\prod_{l=1}^{j}}}(\lam_{i_{l}}\IR-\LR_{i_{l}})\RR_{\lam})f)|\\
\leq &c(|\ti{\lam}_{i_{0}}(\e)|+\|\tLR_{i_{0}}(\e,\cd)\|_{\BR^{i_{0}+1}\to \BR^{0}})\|f\|_{k+1}.
}
for some constant $c>0$.
Hence $\|\ti{\kappa}_{k}(\e,\cd)\|_{\BR^{k+1}_{1}\to \K}$ vanishes as $\e\to 0$.
\proe
Next we consider convergence of the remainder of asymptotic expansion of the eigenfunction under the conditions (I), (II) and (IV). We introduce the notion of weak boundedness. 
\defi
{\label{def:weakbdd}
Let $(\BR^{0},\|\cd\|_{0}),(\BR^{1},\|\cd\|_{1})$ be Banach spaces over $\K$ with $\BR^{1}\subset \BR^{0}$ as linear subspace. We say a linear operator $\LR\,:\,\BR^{1}\to \BR^{0}$ {\it weak bounded} if for any $\e>0$ there exists a constant $c_{\adl{wb}}(\e)>0$ such that
\ali
{
\|\LR f\|_{0}\leq c_{\adr{wb}}(\e)\|f\|_{0}+\e\|f\|_{1}
}
for any $f\in \BR^{1}$. Observe that if $\LR\,:\,(\BR^{1},\|\cd\|_{0})\to (\BR^{0},\|\cd\|_{0})$ is bounded, then this operator is weak bounded.
}
\rem
{
The resolvent $(\lam\IR-\RR)^{-1}\,:\,\BR^{1}\to \BR^{0}$ in \cite[Theorem 8.1]{GL} is weak bounded (see also the proof of Theorem \ref{th:GLtoGAPT}).
}
\prop
{\label{prop:prop_weakbdd}
Let $(\BR^{0},\|\cd\|_{0}),(\BR^{1},\|\cd\|_{1})$ be Banach spaces over $\K$ with $\BR^{1}\subset \BR^{0}$ as linear subspace, and $\LR\,:\,\BR^{1}\to \BR^{0}$ a linear operator.
Then the following (a)-(c) are equivalent:
\ite
{
\item[(a)] The operator $\LR$ is weak bounded.
\item[(b)] Whenever $f_{\e}, f\in \BR^{1}$ satisfy $\|f_{\e}-f\|_{0}\to 0$ as $\e\to 0$ and $\sup_{\e>0}\|f_{\e}\|_{1}<+\infty$, then $\|\LR f_{\e}-\LR f\|_{0}\to 0$ as $\e\to 0$.
\item[(c)] $\sup_{f\in \BR^{1}\,:\,\|f\|_{1}\leq 1,\ \|f\|_{0}\leq \delta}\|\LR f\|_{0}\to 0$ as $\delta\to 0$.
}
}
\pros
Firstly we assume that $\LR$ is weak bounded. Then for any $f_{\e}, f\in \BR^{1}$ satisfying $\|f_{\e}-f\|_{0}\to 0$ as $\e\to 0$ and $c_{\adl{bdfe}}:=\sup_{\e>0}\|f_{\e}\|_{1}<\infty$, we have that for any small $\eta>0$
\ali
{
\|\LR(f_{\e}-f)\|_{0}\leq &c_{\adr{wb}}(\eta)\|f_{\e}-f\|_{0}+\eta\|f_{\e}-f\|_{1}\\
\leq&c_{\adr{wb}}(\eta)\|f_{\e}-f\|_{0}+c_{\adr{bdfe}}\eta+\eta\|f\|_{1}\to c_{\adr{bdfe}}\eta+\eta\|f\|_{1}
}
as $\e\to 0$. 
The arbitrariness of $\eta$ implies the assertion (b).
\smallskip
\par
Secondly we prove that (b) implies (c). To do this, we assume (b) and suppose that (c) does not hold. Then there exists $\e>0$ such that for any $n\geq 1$, there is $f_{n}\in \BR^{1}$ with $\|f_{n}\|_{0}\leq 1/n$ and $\|f_{n}\|_{1}\leq 1$ so that $\|\LR f_{n}\|_{0}\geq \e$. In this case, $\|f_{n}\|_{0}\to 0$ as $n\to \infty$ and therefore $\|\LR f_{n}\|_{0}$ must vanish by the assumption (b). This contradicts with $\|\LR f_{n}\|_{0}\geq \e$. Thus (b) implies (c).
\smallskip
\par
Finally we assume (c). Let $\eta>0$. Then there exists a constant $\delta>0$ such that whenever $g\in \BR^{1}$ satisfies $\|g\|_{0}\leq \delta$ and $\|g\|_{1}\leq 1$, we have $\|\LR g\|_{0}\leq \eta$. Choose any $f\in \BR^{1}$ with $f\neq 0$. Let $c_{\adl{pwd1}}=\|f\|_{0}$ and $c_{\adl{pwd2}}=\|f\|_{1}$. Put $c_{\adl{pwd3}}=\min((\delta c_{\adr{pwd2}})/c_{\adr{pwd1}}, 1)$. In this case, we see
\ali
{
\left\|\frac{c_{\adr{pwd3}}}{c_{\adr{pwd2}}}f\right\|_{0}\leq \delta\quad\text{ and }\quad \left\|\frac{c_{\adr{pwd3}}}{c_{\adr{pwd2}}}f\right\|_{1}\leq 1.
}
Therefore the assumption yields $\|\LR ((c_{\adr{pwd3}}/c_{\adr{pwd2}})f)\|_{0}\leq \eta$ and thus
\ali
{
\|\LR f\|_{0}\leq \eta\frac{c_{\adr{pwd2}}}{c_{\adr{pwd3}}}= &\eta c_{\adr{pwd2}}\max(\frac{c_{\adr{pwd1}}}{\delta c_{\adr{pwd2}}},1)
=\eta\max(\frac{c_{\adr{pwd1}}}{\delta},c_{\adr{pwd2}})
\leq\frac{\eta}{\delta}\|f\|_{0}+\eta\|f\|_{1}.
}
Of course, this inequality holds for $f=0$. Hence $\LR$ is weak bounded.
\proe
\thms
\label{th:asymp_efunc_nB_gene}
Assume that the conditions (I), (II) and (IV) are satisfied. We also assume the following (a)-(c):
\ite
{
\item[(a)] Each $\BR^{k}$ is a Banach space endowed.
\item[(b)] $\nu\in \LR(\BR^{0},\K)$, $\LR_{j}\in \LR(\BR^{i},\BR^{i-j})$, $\sup_{\e>0}\|\tLR_{j}(\e,\cd)\|_{\BR^{i}\to \BR^{i-j}}<+\infty$ for $j=0,1,\dots, n$ and $i=j,j+1,\cdots, n+1$, $\sup_{\e>0}\|g(\e,\cd)\|_{\BR^{i}}<+\infty$ for $i=0,1,\dots, n+1$, $\RR_{\lam}\in \LR(\BR^{1},\BR^{0})\cap \LR(\BR^{i})$ for $i=1,2,\dots, n+1$. Furthermore $\RR_{\lam}\,:\,\BR^{1}\to \BR^{0}$ is weak bounded.
\item[(c)] $\lim_{\e\to 0}\|\tLR_{j}(\e,\cd)\|_{\BR^{i+1}\to\BR^{i-j}}=0$ for each $j=0,1,\dots, n$ and $i=j,\cdots, n$.
}
Then
$\|\ti{g}_{n}(\e,\cd)\|_{\BR^{0}}\to 0$. In particular, $\|\ti{g}_{k}(\e,\cd)\|_{\BR^{i}}\to 0$ for $k,i\geq 0$ with $k+i\leq n$ and $\sup_{\e>0}\|\ti{g}_{k}(\e,\cd)\|_{\BR^{n+1-k}}<\infty$ for $k=0,1,\dots, n$.
\thme
\pros
In order to prove the assertion, we begin with the weak boundedness of $\SR\,:\,\BR^{1}\to \BR^{0}$.
Recall the form $\SR=\RR_{\lam}(\IR-\PR)$. We see the boundedness
\ali
{
\|\SR f\|_{\BR^{0}}\leq (c_{\adr{wb}}(\eta)+\|\RR_{\lam}\|_{\BR^{1}\to \BR^{0}}\|h\|_{\BR^{1}}\|\nu\|_{\BR^{0}\to \K})\|f\|_{\BR^{0}}+\eta \|f\|_{\BR^{1}}
}
for any $f\in \BR^{1}$ and $\eta>0$. Then the operator $\SR\,:\,\BR^{1}\to \BR^{0}$ is weak bounded. Similarity
\alil
{
\|\SR f\|_{\BR^{j}}\leq c_{\adr{dresjj}}c_{\adr{bhi}}c_{\adr{bnui}}\|f\|_{\BR^{0}}+c_{\adr{dresjj}}\|f\|_{\BR^{j}}\label{eq:bSBj0}
}
for $j=1,2,\dots, n+1$ by letting $c_{\adl{dresjj}}=\max_{1\leq j\leq n}\|\RR_{\lam}\|_{\BR^{j}}$, $c_{\adl{bnui}}=\|\nu\|_{\BR^{0}\to \K}$ and $c_{\adl{bhi}}=\max_{0\leq j\leq n+1}\|h\|_{\BR^{j}}$. 
Furthermore we see
\alil
{
\|\TR(\e,\cd)\tLR_{j}(\e,g(\e,\cd))\|_{\BR^{i}}\leq& c_{\adr{bnui}}c_{\adr{dhei}}\|\tLR_{j}(\e,\cd)\|_{\BR^{j+1}\to \BR^{0}}+c_{\adr{dhei}}\|\tLR_{j}(\e,\cd)\|_{\BR^{j+i+1}\to \BR^{i}}=o(1)\label{eq:tL0gB0}\\
\|\TR(\e,\cd)\tLR_{j}(\e,g(\e,\cd))\|_{\BR^{n+1-j}}&\leq c_{\adr{bnui}}c_{\adr{dhei}}\|\tLR_{j}(\e,\cd)\|_{\BR^{j+1}\to \BR^{0}}+c_{\adr{dhei}}\|\tLR_{j}(\e,\cd)\|_{\BR^{n+1}\to \BR^{n+1-j}}\label{eq:tLjgB1}\\
&\leq c_{\adr{bnui}}c_{\adr{dhei}}c_{\adr{ntLje}}+c_{\adr{dhei}}c_{\adr{ntLje}}.\nonumber
}
for $j=0,1,\dots, n$ and $i=0,1,\dots, n-j$, where $c_{\adl{dhei}}=\max_{0\leq i\leq n+1}\|g(\e,\cd)\|_{\BR^{i}}$ and $c_{\adl{ntLje}}=\max_{0\leq j\leq n}\max_{j\leq i\leq n+1}\|\tLR_{j}(\e,\cd)\|_{\BR^{i}\to\BR^{i-j}}$.
Firstly we show the assertion in the case $k=0$.
\ncla
{
$\|\ti{g}_{0}(\e,\cd)\|_{\BR^{i}}\to 0$ for $0\leq i\leq n$ and $\sup_{\e>0}\|\ti{g}_{0}(\e,\cd)\|_{\BR^{n+1}}<+\infty$.
}
Recall the form $\ti{g}_{0}(\e,\cd)=\SR\TR(\e,\cd)\tLR_{0}(\e,g(\e,\cd))$ in (\ref{eq:tgke=...in}) with $k=0$. Therefore, if $n=0$ then (\ref{eq:tL0gB0}) with $j=i=0$ and (\ref{eq:tLjgB1}) with $j=0$ and the weak boundedness of $\SR$ imply $\ti{g}_{0}(\e,\cd)\to 0$ in $\BR^{0}$. If $n\geq 1$ then the boundedness (\ref{eq:bSBj0}) of $\SR$ and (\ref{eq:tL0gB0}) follow $\ti{g}_{0}(\e,\cd)\to 0$ in $\BR^{i}$ for $0\leq i\leq n$. Finally (\ref{eq:bSBj0}) with $j=n+1$ and (\ref{eq:tLjgB1}) with $j=0$ yield the boundedness of $\ti{g}_{0}(\e,\cd)$ in $\BR^{n+1}$ uniformly in $\e>0$.
\cla
{
The assertion holds in the case $n\geq 1$.
}
In order to see this, we assume that $k\geq 1$, $\|\ti{g}_{j}(\e,\cd)\|_{\BR^{i}}\to 0$ and $\|\ti{g}_{j}(\e,\cd)\|_{\BR^{n+1-j}}\leq c_{\adr{btgjBi}}$ for any $0\leq j\leq k-1$ and $0\leq i\leq n-k$ for some $c_{\adl{btgjBi}}>0$. We will check the case $j=k$.
In view of the form (\ref{eq:tgke=...in}) of $\ti{g}_{k}(\e,\cd)$, we obtain the estimate
\ali
{
\|(\TR(\e,\cd)\LR_{l}+\lam_{l}\IR)\ti{g}_{k-l}(\e,\cd)\|_{\BR^{i}}\leq& \|g(\e,\cd)\|_{\BR^{i}}\|\nu\|_{\BR^{0}\to \K}\|\LR_{l}\|_{\BR^{l}\to \BR^{0}}\|\ti{g}_{k-l}(\e,\cd)\|_{\BR^{0}}\\
&+\|\LR_{l}\|_{\BR^{l+i}\to \BR^{i}}\|\ti{g}_{k-l}(\e,\cd)\|_{\BR^{i+l}}+|\lam_{l}|\|\ti{g}_{k-l}(\e,\cd)\|_{\BR^{i}}
}
for $k,i\geq 0, 1\leq l\leq k$ with $k+i\leq n$. It follows from the assumption that $(\TR(\e,\cd)\LR_{l}+\lam_{l}\IR)\ti{g}_{k-l}(\e,\cd)$ vanishes in $\BR^{i}$ in the case $k+i\leq n$. Similarity, this function is bounded in $\BR^{n+1-k}$. Thus together with (\ref{eq:tL0gB0}) and (\ref{eq:tLjgB1}), we obtain convergence $\|\ti{g}_{k}(\e,\cd)\|_{\BR^{i}}\to 0$ for $0\leq i\leq n-k$ and the boundedness of $\|\ti{g}_{k}(\e,\cd)\|_{\BR^{n+1-k}}$. 
Hence we obtain the assertion inductively.
\proe
\rem
{
Theorem \ref{th:asymp_efunc_nB_gene} is a generalization of \cite[Theorem 2.4]{T2011} to general operators.
}
Finally we consider an asymptotic behaviour of the eigenfunction with sequence of Banach spaces that depend on $\e>0$.
\thms
\label{th:asymp_efunc_nospecgap}
Let $\BR^{0}_{\e},\BR^{1}_{\e},\cdots, \BR^{n+1}_{\e}$ be a sequence of Banach spaces over $\K$ with $\BR^{0}_{\e}\supset \cdots \supset \BR^{n+1}_{\e}$ for each $\e>0$. We assume the following (a)-(c):
\ite
{
\item[(a)] Operators $\LR,\LR_{1},\dots, \LR_{n},\LR(\e,\cd)$ satisfy the conditions (I)-(IV) replacing as $\BR^{k}=\BR^{k}_{\e}$ $(k=0,1,\dots,n+1)$ for each $\e>0$.
\item[(b)] There exist constants $c_{\adl{bhBj}},c_{\adl{bnBj}},c_{\adl{bLjBj}},c_{\adl{bnueBj}}>0$ such that $\|h\|_{\BR^{j}_{\e}}\leq c_{\adr{bhBj}}$, $\|\nu\|_{\BR^{j}_{\e}\to \K}\leq c_{\adr{bnBj}}$, $\|\LR_{j}\|_{\BR^{i}_{\e}\to \BR^{i-j}_{\e}}\leq c_{\adr{bLjBj}}$ and $\|\kappa(\e,\cd)\|_{\BR^{j}_{\e}\to \K}\leq c_{\adr{bnueBj}}$ for all $\e>0$, $j=0,1,\dots, n$ and $i=j,\cdots, n+1$.
\item[(c)] For any $0\leq k\leq n$, $2\leq p\leq n+2$ and $0\leq k_{1}\leq \cdots \leq k_{p}\leq k$, 
\ali
{
\textstyle{(\prod_{q=1}^{p-1}\|\SR\|_{\BR^{k_{q}}_{\e}})}\|\tLR_{k_{p}-k_{p-1}}(\e,\cd)\|_{\BR^{1+k_{p}}_{\e}\to \BR^{k_{p-1}}_{\e}}\to 0
}
as $\e\to 0$ for each $0\leq k\leq n$.
}
Then $\|\ti{g}_{k}(\e,\cd)\|_{\BR^{i}_{\e}}$ vanishes as $\e\to 0$ for any $0\leq k\leq n$ and $0\leq i\leq n-k$.
\thme
\pros
\ncla
{
If $\dim \BR^{n+1}_{\e}=1$, then $g(\e,\cd)$ equals $h$, namely the assertion holds.
}
Indeed, since $h$ is in $\BR^{n+1}_{\e}$, its space forms $\BR^{n+1}_{\e}=\{ch\,:\,c\in \K\}$. The function $h(\e,\cd)$ belongs also to $\BR^{n+1}_{\e}$ and therefore $h(\e,\cd)=c(\e)h$ for some constant $c(\e)\neq 0$. Thus $g(\e,\cd)=h(\e,\cd)/\nu(h(\e,\cd))$ is equal to $h$.
\smallskip
\par
In what follows, we assume $\dim \BR^{n+1}_{\e}\geq 2$. 
\cla
{
For some $\e_{0}>0$, the number $c_{\adl{bdge}}:=\max_{0\leq i\leq n}\sup_{0<\e<\e_{0}}\|g(\e,\cd)\|_{\BR^{i}_{\e}}$ is finite.
}
Indeed, it follows from the assumption (b) that
\ali
{
|\ti{\lam}_{0}(\e)|&=|\kappa(\e,\tLR_{0}(\e,h))|
\leq c_{\adr{bhBj}}c_{\adr{bnueBj}}\|\tLR_{0}(\e,\cd)\|_{\BR^{i}_{\e}}
}
Therefore, the form (\ref{eq:tgke=...in2}) in $k=0$ of $\ti{g}_{0}(\e,\cd)$ implies
\ali
{
|\|g(\e,\cd)\|_{\BR^{i}_{\e}}-\|h\|_{\BR^{i}_{\e}}|\leq\|\ti{g}_{0}(\e,\cd)\|_{\BR^{i}_{\e}}=&\|-\SR\tLR_{0}(\e,g(\e,\cd))+\ti{\lam}_{0}(\e)\SR g(\e,\cd)\|_{\BR^{i}_{\e}}\\
\leq&(1+c_{\adr{bhBj}}c_{\adr{bnueBj}})\|\SR\|_{\BR^{i}_{\e}}\|\tLR_{0}(\e,\cd)\|_{\BR^{i}_{\e}}\|g(\e,\cd)\|_{\BR^{i}_{\e}}.
}
Thus
\alil
{
|1-\|h\|_{\BR^{i}_{\e}}/\|g(\e,\cd)\|_{\BR^{i}_{\e}}|\leq(1+c_{\adr{bhBj}}c_{\adr{bnueBj}})\|\SR\|_{\BR^{i}_{\e}}\|\tLR_{0}(\e,\cd)\|_{\BR^{i}_{\e}}\label{eq:ineq_hge}
}
vanishes as $\e\to 0$ by the assumption (c). Choose any $1/2>\eta>0$, there exists $\e_{0}>0$ such that for any $0<\e<\e_{0}$, $-1/2<-\eta<1-\|h\|_{\BR^{i}_{\e}}/\|g(\e,\cd)\|_{\BR^{i}_{\e}}<\eta<1/2$ and so $(3/2)^{-1}<\|g(\e,\cd)\|_{\BR^{i}_{\e}}/\|h\|_{\BR^{i}_{\e}}<(1/2)^{-1}$.
We see
\ali
{
\|g(\e,\cd)\|_{\BR^{i}_{\e}}=\|h\|_{\BR^{i}_{\e}}\frac{\|g(\e,\cd)\|_{\BR^{i}_{\e}}}{\|h\|_{\BR^{i}_{\e}}}< 2\|h\|_{\BR^{i}_{\e}}\leq 2\sup_{0<\e<\e_{0}}\|h\|_{\BR^{i}_{\e}}.
}
Hence the claim is satisfied.
\cla
{\label{cla:pos_S}
$\inf_{1\leq i\leq n+1,\e>0}\|\SR\|_{\BR^{i}_{\e}}>0$.
}
Note that the dimension of $\BR^{i}_{\e}$ is larger than $1$ for any $1\leq i\leq n+1$. Therefore there is $f\in \BR^{i}_{\e}$ such that $g:=(\IR-\PR)f\neq 0$. Then $\|(\IR-\PR)g\|_{\BR^{i}_{\e}}=\|g\|_{\BR^{i}_{\e}}$ and thus $1\leq \|(\IR-\PR)\|_{\BR^{i}_{\e}}$. Note that $\IR-\PR=(\RR-\lam\IR)\SR$ is an operator from $\BR^{i}_{\e}$ to $\BR^{i}_{\e}$. We see
\ali
{
1\leq \|\IR-\PR\|_{\BR^{i}_{\e}}&\leq \|\RR-\lam\IR\|_{\BR^{i}_{\e}}\|\SR\|_{\BR^{i}_{\e}}\\
&=\|\LR-\lam\PR-\lam\IR\|_{\BR^{i}_{\e}}\|\SR\|_{\BR^{i}_{\e}}\\
&\leq (c_{\adr{bLjBj}}+|\lam|(c_{\adr{bhBj}}c_{\adr{bnBj}}+1))\|\SR\|_{\BR^{i}_{\e}}
}
and thus $\|\SR\|_{\BR^{i}_{\e}}\geq (c_{\adr{bLjBj}}+|\lam|(c_{\adr{bhBj}}c_{\adr{bnBj}}+1))^{-1}$.
\cla
{
The assertion holds.
}
Assume $0<\e<\e_{0}$. We recall the form (\ref{eq:tgke=...in2}) of $\ti{g}_{k}(\e,\cd)$.
By the form $\SR\TR(\e,f)=\nu(f)\SR g(\e,\cd)-\SR f$, we note the estimate 
\ali
{
\|\SR\TR(\e,f)\|_{\BR^{i}_{\e}}\leq& |\nu(f)|\|\SR\|_{\BR^{i}_{\e}}\|g(\e,\cd)\|_{\BR^{i}_{\e}}+\|\SR f\|_{\BR^{i}_{\e}}\\
\leq&\|\SR\|_{\BR^{i}_{\e}}(c_{\adr{bdge}}\|\nu\|_{\BR^{i}\to\K}+1)\|f\|_{\BR^{i}_{\e}}
\leq c_{\adr{sect}}\|\SR\|_{\BR^{i}_{\e}}\|f\|_{\BR^{i}_{\e}}
}
by putting $c_{\adl{sect}}=c_{\adr{bdge}}c_{\adr{bnBj}}+1$.
Let $\MR^{0}_{i}=\lam_{i}\SR$ and $\MR^{1}_{i}=\SR\TR(\e,\cd)\LR_{i}$. We rewrite $\ti{g}_{k}(\e,\cd)$ by
\ali
{
\ti{g}_{k}(\e,\cd)=&\SR\TR(\e,\cd)\tLR_{k}(\e,g(\e,\cd))\\
&+\sum_{p=2}^{k+1}\sum_{\underset{i_{1}+\cdots+i_{p}=k}{i_{1},\dots,i_{p-1}\geq 1,i_{p}\geq 0\,:}}\sum_{j_{1},\dots, j_{p-1}\in \{0,1\}}\MR^{j_{1}}_{i_{1}}\MR^{j_{2}}_{i_{2}}\cdots\MR^{j_{p-1}}_{i_{p-1}}\SR\TR(\e,\cd)\tLR_{i_{p}}(\e,g(\e,\cd)).
}
Note also that 
\ali
{
\|\MR^{0}_{i}f\|_{\BR^{k}_{\e}}\leq& |\lam_{i}|\|\SR\|_{\BR^{k}_{\e}}\|f\|_{\BR^{k}_{\e}},\qquad
\|\MR^{1}_{i}f\|_{\BR^{k}_{\e}}\leq c_{\adr{bLjBj}}c_{\adr{sect}}\|\SR\|_{\BR^{k}_{\e}}\|f\|_{\BR^{k+i}_{\e}}.
}
Namely, $\|\MR^{j}_{i}\|_{\BR^{k}_{\e}\to\BR^{k-i\cd j}_{\e}}\leq c_{\adl{bMji}}\|\SR\|_{\BR^{k}_{\e}}$ for $j=0,1$ and $0\leq i\leq n$ with $k\geq ij$ by putting $c_{\adr{bMji}}=\max\{\max_{0\leq i\leq n}|\lam_{i}|,c_{\adr{bLjBj}}c_{\adr{sect}}\}$. Consequently, we obtain that for $0\leq i\leq n-k$
\ali
{
\|\SR\TR(\e,\cd)\tLR_{k}(\e,g(\e,\cd))\|_{\BR^{i}_{\e}}\leq c_{\adr{sect}}c_{\adr{bdge}}\|\SR\|_{\BR^{i}_{\e}}\|\tLR_{k}(\e,\cd)\|_{\BR^{k+i+1}_{\e}\to\BR^{i}_{\e}}\to 0
}
and by putting $k_{p}=j_{1}i_{1}+\cdots+j_{p}i_{p}$
\ali
{
&\|\MR^{j_{1}}_{i_{1}}\MR^{j_{2}}_{i_{2}}\cdots\MR^{j_{p-1}}_{i_{p-1}}\SR\TR(\e,\cd)\tLR_{i_{p}}(\e,g(\e,\cd))\|_{\BR^{i}_{\e}}\\
\leq&c_{\adr{bdge}}\|\MR^{j_{1}}_{i_{1}}\|_{\BR^{i+k_{1}}_{\e}\to\BR^{i}_{\e}}\Big(\prod_{l=1}^{p-1}\|\MR^{j_{2}}_{i_{2}}\|_{\BR^{i+k_{l+1}}\to\BR^{i+k_{l}}}\Big)\|\tLR_{i_{p}}(\e,\cd)\|_{\BR^{1+i+k_{p-1}+i_{p}}_{\e}\to\BR^{i+k_{p-1}}_{\e}}\\
\leq&c_{\adr{bMji}}^{p}c_{\adr{bdge}}\|\SR\|_{\BR^{i}_{\e}}\|\SR\|_{\BR^{i+k_{1}}_{\e}}\|\SR\|_{\BR^{i+k_{2}}_{\e}}\cdots\|\SR\|_{\BR^{i+k_{p-1}}_{\e}}\|\tLR_{i_{p}}(\e,\cd)\|_{\BR^{1+i+k_{p-1}+i_{p}}_{\e}\to\BR^{i+k_{p-1}}_{\e}}\to 0
}
as $\e\to 0$ by (c). Hence the proof is complete.
\proe
\section{Examples}\label{sec:ex}
\subsection{Coefficients and remainders of asymptotic expansions in the case $n=0,1,2$}\label{sec:coefn=012}
Referring to the notation given in Section \ref{sec:rep_remainer}, we have
\ali
{
\lam_{1}=&\nu(\LR_{1}h),\ \lam_{2}=\nu(\LR_{2}h)+\nu((\lam_{1}\IR-\LR_{1})\RR_{\lam}\LR_{1}h)\\
\kappa_{1}=&((\lam_{1}\IR-\LR_{1})\RR_{\lam})^{*}\nu,\ \ \kappa_{2}=((\lam_{1}\IR-\LR_{1})\RR_{\lam}(\lam_{1}\IR-\LR_{1})\RR_{\lam})^{*}\nu+((\lam_{2}\IR-\LR_{2})\RR_{\lam})^{*}\nu\\
g_{1}=&\SR(\lam_{1}\IR-\LR_{1})h,\ \ g_{2}=\SR(\lam_{2}\IR-\LR_{2})h+\SR(\lam_{1}\IR-\LR_{1})\SR(\lam_{1}\IR-\LR_{1})h\\
\ti{\lam}_{0}(\e)=&\kappa(\e,\tLR_{0}(\e,h)),\ \ \ti{\lam}_{1}(\e)=\kappa(\e,\tLR_{1}(\e,h)+\tLR_{0}(\e,\cd)\QR(\e,\cd)\RR_{\lam}\LR_{1}h)\\
\ti{\lam}_{2}(\e)=&\kappa(\e,\tLR_{2}(\e,h)+\tLR_{0}(\e,\cd)\QR(\e,\cd)\RR_{\lam}\LR_{2}h)+\\
&+\kappa(\e,\tLR_{1}(\e,\cd)\QR(\e,\cd)\RR_{\lam}\LR_{1}h)+\kappa(\e,\tLR_{0}(\e,\cd)\QR(\e,\cd)\RR_{\lam}(\LR_{1}\QR(\e,\cd)\RR_{\lam}+\lam_{1}\RR_{\lam})\LR_{1}h)\\
\ti{\kappa}_{0}(\e,\cd)=&(\tLR_{0}(\e,\cd)\QR(\e,\cd)\RR_{\lam})^{*}\kappa(\e,\cd)\\
\ti{\kappa}_{1}(\e,\cd)=&(\tLR_{1}(\e,\cd)\QR(\e,\cd)\RR_{\lam}+\tLR_{0}(\e,\cd)\QR(\e,\cd)\RR_{\lam}(\LR_{1}\QR(\e,\cd)\RR_{\lam}+\lam_{1}\RR_{\lam}))^{*}\kappa(\e,\cd)\\
\ti{\kappa}_{2}(\e,\cd)=&(\tLR_{2}(\e,\cd)\QR(\e,\cd)\RR_{\lam}+\sum_{l=0}^{1}\tLR_{l}(\e,\cd)\QR(\e,\cd)\RR_{\lam}(\LR_{2-l}\QR(\e,\cd)\RR_{\lam}+\lam_{2-l}\RR_{\lam}))^{*}\kappa(\e,\cd)\\
&+(\tLR_{0}(\e,\cd)\QR(\e,\cd)\RR_{\lam}(\LR_{1}\QR(\e,\cd)\RR_{\lam}+\lam_{1}\RR_{\lam})(\LR_{1}\QR(\e,\cd)\RR_{\lam}+\lam_{1}\RR_{\lam}))^{*}\kappa(\e,\cd)\\
\ti{g}_{0}(\e,\cd)=&\SR\TR(\e,\cd)\tLR_{0}(\e,g(\e,\cd))\\
\ti{g}_{1}(\e,\cd)=&\SR\TR(\e,\cd)\tLR_{1}(\e,g(\e,\cd))+(\SR\TR(\e,\cd)\LR_{1}+\lam_{1}\SR)\SR\TR(\e,\cd)\tLR_{0}(\e,g(\e,\cd))\\
\ti{g}_{2}(\e,\cd)=&\SR\TR(\e,\cd)\tLR_{2}(\e,g(\e,\cd))+
\sum_{l=0}^{1}(\SR\TR(\e,\cd)\LR_{2-l}+\lam_{2-l}\SR)\SR\TR(\e,\cd)\tLR_{l}(\e,g(\e,\cd))\\
&+(\SR\TR(\e,\cd)\LR_{1}+\lam_{1}\SR)(\SR\TR(\e,\cd)\LR_{1}+\lam_{1}\SR)\SR\TR(\e,\cd)\tLR_{0}(\e,g(\e,\cd)).
}
\subsection{Asymptotic behaviour of Ruelle operators for countable Markov shifts}\label{ex:CMS}
The aim of this section is to extend the asymptotic perturbation theory of Ruelle operators in the finite state case \cite{T2011} to the infinite state case.
We use the notation in Section \ref{sec:rep_remainer} and the notion of topological Markov shifts and Ruelle operators given in Section \ref{sec:thermo}. 
\smallskip
\par
Let $X$ be a topological Markov shift with countable state space $S$ and with finitely irreducible transition matrix. Take functions $\ph, \ph_{1},\ph_{2},\dots,\ph_{n}\in F_{\theta}(X,\R)$ with $P(\ph)<\infty$ and $\theta\in (0,1)$, and $\phe\in F_{\theta(\e)}(X,\R)$ with $P(\phe)<\infty$ and $\theta(\e)\in [\theta,1)$ $(0<\e<1)$. 
We put $\ti{\ph}_{k}(\e,\cd):=(\phe-\ph-\ph_{1}\e-\cdots-\ph_{k}\e^{k})/\e^{k}$ for $k=0,1,\dots, n$. Then we have the form
\alil
{
\ph(\e,\cd)=\ph+\ph_{1}\e+\cdots+\ph_{n}\e^{n}+\ti{\ph}_{n}(\e,\cd)\e^{n}.\label{eq:phe=ph+...}
}
By applying Taylor Theorem to the map $\e \mapsto \exp(\phe)$, we obtain the expansion of the Ruelle operator $\LR(\e,\cd)$ of $\phe$
\ali
{
\LR(\e,\cd)=&\LR+\LR_{1}\e+\cdots+\LR_{n}\e^{n}+\tLR_{n}(\e,\cd)\e^{n}
}
by putting
\alil
{
\LR_{k}f(\om)&:=\LR(F_{k}f)(\om)\quad \text{ with}\label{eq:LRk=}\\
F_{k}&:
=\sum_{l_{1},\dots,l_{k}\geq 0\,:\,\atop{l_{1}+2l_{2}+\cdots+kl_{k}=k}}\frac{(\ph_{1})^{l_{1}}\cdots (\ph_{k})^{l_{k}}}{l_{1}!\cdots l_{k}!}\label{eq:Fk=}\\
\tLR_{n}(\e,f)(\om)&:=\LR(\ti{F}_{n}(\e,\cd)f)(\om)\quad\text{ with}\label{eq:tLRke=}\\
\ti{F}_{n}(\e,\cd)&:=
\case
{
e^{\ti{\ph}_{0}(\e,\cd)}-1,&n=0\\
\di\ti{\ph}_{1}(\e,\cd)+\int_{0}^{1}(e^{t\ti{\ph}_{0}(\e,\cd)}-1)\,dt(\ph_{1}+\ti{\ph}_{1}(\e,\cd)),&n=1\\
\di\sum_{l=1}^{n}\frac{1}{l!}\Big\{\sum_{i_{1},\dots,i_{l-1}\geq 1,i_{l}\geq 0\,:\,\atop{i_{1}+\dots+i_{l}=n}}\ph_{i_{1}}\cdots\ph_{i_{l-1}}\ti{\ph}_{i_{l}}(\e,\cd)\\
\di\qquad +\sum_{i=1}^{n}\ph_{i}\ti{\ph}_{n-i}(\e,\cd)\ti{\ph}_{0}(\e,\cd)^{l-2}+\ti{\ph}_{n}(\e,\cd)\ti{\ph}_{0}(\e,\cd)^{l-1}\Big\}\\
\di\qquad +\int_{0}^{1}\frac{(1-t)^{n-1}}{(n-1)!}(e^{t\ti{\ph}_{0}(\e,\cd)}-1)\,dt(\ph_{1}+\ti{\ph}_{1}(\e,\cd))^{n},&n\geq 2.\\
}\label{eq:tFke=}
}
formally for $f\,:\,X\to \C$ and $\om\in X$ for each $k=1,2,\dots, n$, where $\LR$ is the Ruelle operator of the potential $\ph$. We introduce the following condition:
\ite
{
\item[($\Ph$)] The operators $\LR_{k}$ $(k=1,2,\dots, n)$ given in (\ref{eq:LRk=}) are bounded linear operators acting both on $C_{b}(X)$ and on $F_{\theta,b}(X)$. 
}
\rem
{
For example, if $\sharp S<+\infty$, then this condition ($\Ph$) holds automatically. As other example,  the condition is satisfied for physical potentials given in perturbed graph directed Markov system with strongly regular (see (\ref{eq:phe=})).
}
Referring to Theorem \ref{th:Spec_gap}, we take the spectral triplet $(\lam,h,\nu)$ of the Ruelle operator $\LR_{\ph}$ of $\ph$ and the triplet $(\lam(\e),h(\e,\cd),\nu(\e,\cd))$ of the Ruelle operator $\LR_{\phe}$ of $\phe$. Since $\lam$ is in the resolvent set of $\RR:=\LR_{\ph}-\lam h\otimes \nu$, $(\RR-\lam \IR)^{-1}$ is a bounded linear operator acting $F_{\theta,b}(X)$. When we set
\ali
{
\BR^{0}=F_{\theta(\e),b}(X),\ \ \DR^{0}=\BR^{1}=\cdots =\BR^{n+1}=F_{\theta,b}(X)
}
for each $\e>0$, the conditions (I)-(VI) are fulfilled with $1_{\BR}:=1$. Therefore $\lam(\e),h(\e,\cd),\nu(\e,\cd)$ have the asymptotic expansions (\ref{eq:le=...}), (\ref{eq:nue=...}), (\ref{eq:he=...}), respectively. Consequently, the topological pressure $P(\phe)=\log\lam(\e)$ and the Gibbs measure $\mu(\e,\cd)=h(\e,\cd)\nu(\e,\cd)$ for $\phe$ satisfy the following expansions:
\props
\label{th:asymp_evec_infRuelle}
Let $X$ be a topological Markov shift whose transition matrix is finitely irreducible. Assume that $\ph,\ph_{1},\dots, \ph_{n}$ satisfy the condition ($\Ph$). Then
\ali
{
P(\varphi(\epsilon,\cdot))=&P(\varphi)+p_{1}\epsilon+\cdots+p_{n}\epsilon^{n}+\ti{p}_{n}(\e)\e^{n}\\
\mu(\e,f)=&\mu(f)+\mu_{1}(f)\e+\cdots+\mu_{n}(f)\e^{n}+\ti{\mu}_{n}(\e,f)\e^{n}\qquad \text{ for }f\in F_{\theta,b}(X)
}
with
\alil
{
p_{k}=&\sum_{l=1}^{k}\frac{(-1)^{l-1}}{l\cd \lam^{l}}\sum_{i_{1},\dots,i_{l}\geq 1\,:\,\atop{i_{1}+\cdots+i_{l}=k}}\lam_{i_{1}}\cdots \lam_{i_{l}}\qquad (k=1,\dots, n)\\
\mu_{k}(f)=&\sum_{i=0}^{k}\nu_{i}(h_{k-i}f)\\
\ti{p}_{n}(\e)=&\sum_{k=1}^{n}\frac{(-1)^{k-1}}{k\lam^{k}}\sum_{i_{1},\dots,i_{k}\geq 0\,:\,\atop{i_{1}+\cdots+i_{k}=n}}\lam_{i_{1}}\cdots\lam_{i_{k-1}}\ti{\lam}_{i_{k}}(\e)\\
&+\int_{0}^{1}\frac{(t-1)^{n-1}}{(n-1)!}\Big((\lam+t\ti{\lam}_{0}(\e))^{-n}-\lam^{-n}\Big)\,dt\,\Big(\frac{\ti{\lam}_{0}(\e)}{\e}\Big)^{n}\nonumber\\
\ti{\mu}_{n}(\e,f)=&\sum_{k=0}^{n}\nu_{k}(\ti{h}_{n-k}(\e,\cd)f)+\ti{\nu}_{n}(\e,h(\e,\cd)f).\label{eq:tnumef=...}
}
\prope
\pros
The former expansion is due to Taylor Theorem for $x\mapsto \log x$ and the later follows directly from the expansions (\ref{eq:nue=...}) and (\ref{eq:he=...}).
\proe
By using convergence conditions given in Theorem \ref{th:asymp_eval_nB_gene2}, Theorem \ref{th:asymp_eval_nB_gene3}, Theorem \ref{th:asymp_efunc_nB_gene} and Theorem \ref{th:asymp_efunc_nospecgap}, we obtain the following sufficient conditions for asymptotic behaviours of $P(\phe)$ and $\mu(\e,\cd)$:
\thms
\label{th:asymp_pphe}
Let $X$ be a topological Markov shift whose transition matrix is finitely irreducible. Assume that $\ph,\ph_{1},\dots, \ph_{n}$ satisfy the condition ($\Ph$). Assume also $\|\tLR_{n}(\e,f)\|_{C}\to 0$ for each $f\in F_{\theta,b}(X)$. Then $\ti{\lam}_{n}(\e)\to 0$ and $\ti{\nu}_{n}(\e,f)\to 0$ for each $f\in \DR^{n}$, where $\DR^{n}$ is defined inductively by (\ref{eq:Dk=...}) and we put $\DR^{0}:=(\RR-\lam\IR)\DR$ with the completion $\DR$ of $F_{\theta,b}(X)$ with respect to the norm $\|\cd\|_{C}$. In particular, $\ti{p}_{n}(\e)\to 0$.
\thme
\pros
Note that $\LR(\e,\cd)$, $\LR_{\ph}$, $\PR:=h\otimes \nu$ and $\RR:=\LR_{\ph}-\lam\PR$ are bounded linear operators acting on $C_{b}(X)$.
Let
\ali
{
\BR^{0}=C_{b}(X),\quad \BR^{1}=\cdots=\BR^{n+1}=F_{\theta,b}(X).
}
 
Since $\lam$ is in the resolvent set of $\RR\in \LR(F_{\theta,b}(X))$, the inclusion $\BR^{1}\subset \DR^{0}$ is valid.
The conditions (I) and (III) are clearly satisfied. Now we check the condition (II). Then we show the invectively of $\RR-\lam\IR$ on $\DR$. In order to show this, we recall the quasi-compactness of $\LR_{\ph}$ acting on $F_{\theta,b}(X)$ from Theorem \ref{th:Spec_gap}: $\LR_{\ph}$ has the spectral decomposition
$\LR_{\ph}=\sum_{i=0}^{p-1}\lam_{i}\PR_{i}+\ER$
satisfying that (a) $\lam_{i}$ is a complex number with $\lam_{0}=\lam$, $\lam_{i}^{p}=\lam^{p}$ and $\lam_{i}\neq\lam_{j}$ for $j\neq i$, (b) $\PR_{i}$ is the projection onto the one-dimensional eigenspace of $\lam_{i}$, and (c) the spectral radius of $\ER$ is less than $\lam$. We notice that $\PR_{i}$ and $\ER$ are extended as bounded linear operators from $\DR$ to $C_{b}(X)$ by $\PR_{i}f=\lim_{m\to \infty}\PR_{i}f_{m}$ and $\ER f=\lim_{m\to \infty}\ER f_{m}$ for any $f_{m}\in F_{\theta,b}(X)$ with $\|f_{m}-f\|_{C}\to 0$. Then these preserve the equations $\PR_{i}\PR_{j}=\PR_{i}\ER=\ER\PR_{i}=\OR$ for $i\neq j$. 
By using the definitions of $\ER$ and the fact that the spectral radius of $\ER\in \LR(F_{\theta,b}(X))$ is less than $\lam$, we obtain
\alil
{
\|\lam^{-n}\ER^{n}f\|_{C}\to 0\label{eq:conv_ER}
}
as $n\to \infty$ for each $f\in \DR$. Now we assume that $(\RR-\lam\IR)f=0$ for an element $f\in \DR$. Since the form $\RR=\sum_{i=1}^{p-1}\lam_{i}\PR_{i}+\ER$ on $\DR$, $\PR_{i}f=\lam^{-1}\PR_{i}\RR f=\lam^{-1}\lam_{i}\PR_{i}f$ namely $\PR_{i}f=0$ for $i\neq 0$. Thus $\ER f=\lam f$. This means $f=0$ by (\ref{eq:conv_ER}). Therefore $\RR-\lam\IR$ is injective on $\DR$ and so the condition (II) is fulfilled.
Since $\|\kappa(\e,\cd)\|_{C\to \K}\leq 1/\inf h$,   $\|\kappa(\e,\tLR_{j}(\e,f))|\leq (1/\inf h)\|\tLR_{j}(\e,f)\|_{C}\to 0$ as $f\in \BR^{0}$, and $|\nu(\e,h)|\leq \|h\|_{C}$ hold, the conditions of Theorem \ref{th:asymp_eval_nB_gene2} are guaranteed. Hence the proof is complete.
\proe
\rem
{
Assume that $X$ is a topological Markov shift whose transition matrix is finitely irreducible. In the case that the state space $S$ is finite, $\|\tLR_{n}(\e,f)\|_{C}\to 0$ for each $f\in F_{\theta,b}(X)$ if and only if $\|\ti{\ph}_{n}(\e,\cd)\|_{C}\to 0$. In the case $\sharp S=+\infty$, there is $\phe$ so that $\|\tLR_{n}(\e,f)\|_{C}$ vanishes for each $f$ even if 
$\|\ti{\ph}_{n}(\e,\cd)\|_{C}\equiv \infty$ (see Section \ref{sec:GDMS_Gibbs}). 
}
\thms
\label{th:asymp_pphe_he}
Let $X$ be a topological Markov shift whose transition matrix is finitely irreducible. We assume that $\ph,\ph_{1},\dots, \ph_{n}$ satisfy the condition ($\Ph$). Assume also that $\|\tLR_{n}(\e,\cd)\|_{F_{\theta}\to C}\to 0$, $\sup_{\e>0}[\phe]_{\theta}^{2}<\infty$ and $\sup_{\e>0}\|\tLR_{n}(\e,\cd)\|_{F_{\theta}}<\infty$. Then we have $\|\ti{h}_{n}(\e,\cd)\|_{L^{1}(\nu)}\to 0$ and $\ti{\mu}_{n}(\e,f)\to 0$ for each $f\in F_{\theta,b}(X)$.
\thme
\pros
Note that $\LR(\e,\cd)$, $\LR_{\ph}$, $\PR:=h\otimes \nu$ and $\RR:=\LR_{\ph}-\lam\PR$ are bounded linear operators acting on $L^{1}(\nu)$. We begin with the following.
\ncla
{
$\|\ti{g}_{n}(\e,\cd)\|_{L^{1}(\nu)}\to 0$ and $\sup_{\e>0}\|\ti{g}_{n}(\e,\cd)\|_{F_{\theta}}<+\infty$.
}
Indeed, we put
\ali
{
\BR^{0}=L^{1}(\nu),\quad \DR^{0}=\BR^{1}=\cdots=\BR^{n+1}=F_{\theta,b}(X).
}
We will see the validity of the conditions of Theorem \ref{th:asymp_efunc_nB_gene}. The conditions (a) in this theorem is clearly satisfied. We will prove that the condition (b) is satisfied. It is clear that $\nu\,:\,\LR^{1}(\nu)\to \R$ is bounded. By using the relations $\|\cd\|_{L^{1}(\nu)}\leq \|\cd\|_{C}\leq \|\cd\|_{F_{\theta}}$, we have
\ali
{
\|\LR_{0}\|_{L^{1}(\nu)}=1<+\infty,\ \|\LR_{j}\|_{L^{1}(\nu)\to F_{\theta}}\leq \|\LR_{j}\|_{F_{\theta}}<+\infty.
}
The boundedness of $\tLR_{j}(\e,\cd)$ follows from
\alil
{
\|\tLR_{n}(\e,\cd)\|_{F_{\theta}\to L^{1}(\nu)}\leq& \|\tLR_{n}(\e,\cd)\|_{F_{\theta}}<+\infty\\
\|\tLR_{k}(\e,\cd)\|_{F_{\theta}\to L^{1}(\nu)}\leq&\|\tLR_{k}\|_{F_{\theta}}\leq \sum_{i=1}^{n-k}\|\LR_{k+i}\|_{F_{\theta}}\e^{i}+\|\tLR_{n}(\e,\cd)\|_{F_{\theta}}\e^{n-k}<+\infty\label{eq:tLke->0}
}
for $k=0,1,\dots, n-1$. To see $\sup_{\e>0}\|g(\e,\cd)\|_{F_{\theta}}<\infty$, we note the form $g(\e,\cd)=h(\e,\cd)/\nu(h(\e,\cd))=f(\e,\cd)/\nu(f(\e,\cd))$ by letting $f(\e,\cd)=h(\e,\cd)/\|h(\e,\cd)\|_{C}$. The inequality
\ali
{
\var_{k}(\log h(\e,\cd))\leq&\sum_{i=k+1}^{\infty}\var_{i}\ph(\e,\cd)
\leq\sum_{i=k+1}^{\infty}\theta^{i}[\ph(\e,\cd)]_{\theta}^{2}
\leq c_{\adr{bde.func}}\theta^{k}
}
is satisfied for each $k\geq 1$ (see \cite[p.67]{Sar09}) with $c_{\adl{bde.func}}=\theta\sup_{\e>0}[\ph(\e,\cd)]_{\theta}^{2}/(1-\theta)$. Therefore $\sup_{\e>0}[\log h(\e,\cd)]_{\theta}^{1}<\infty$. Moreover, we have for $\om,\up\in X$ with $\om_{0}=\up_{0}$, 
\ali
{
|h(\e,\om)-h(\e,\up)|=&h(\e,\up)|e^{\log h(\e,\om)-\log h(\e,\up)}-1|\\
=&h(\e,\up)e^{\alpha(\log h(\e,\om)-\log h(\e,\up))}|\log h(\e,\om)-\log h(\e,\up)|\\
\leq&h(\e,\up)e^{\theta c_{\adr{bde.func}}}c_{\adr{bde.func}}d_{\theta}(\om,\up).
}
Thus $\{f(\e,\cd)\,:\,\e>0\}$ is equicontinuous and uniformly bounded. Choose a positive sequence $(\e_{n})$ with $\lim_{n\to \infty}\e_{n}=0$ so that $\lim_{n\to \infty}\|g(\e_{n},\cd)\|_{C}=\limsup_{\e\to 0}\|g(\e,\cd)\|_{C}$. Ascoli Theorem implies that there exist a subsequence $(\e^\p_{n})$ of $(\e_{n})$ and a continuous function $f\,:\,X\to \R$ such that $f(\e_{n}^\p,\cd)$ converges to $f$ for each point of $X$ as $n\to \infty$. By the Lebesgue dominated convergence theorem, 
the inequality (\ref{eq:ineq_hge}) implies $f(\e_{n}^\p,\cd)$ converges to $f$ in $L^{1}(X,\BR(X),\nu)$. Since $f$ is continuous (in particular, in $F_{\theta,b}(X)$), $\|f\|_{C}=1$ implies $\nu(f)>0$. This means that $g(\e,\cd)=f(\e,\cd)/\nu(f(\e,\cd))$ converges to $f/\nu(f)$ in $L^{1}(X,\BR(X),\nu)$ running through $\e\in \{\e_{n}^\p\}$. We notice $\|g(\e,\cd)\|_{C}\leq 1/\nu(f(\e,\cd))$ and then 
\ali
{
\limsup_{\e\to 0}\|g(\e,\cd)\|_{C}=\lim_{n\to \infty}\|g(\e_{n}^\p,\cd)\|_{C}\leq 1/\nu(f)
}
Therefore, $c_{\adl{bdge2}}=\sup_{0<\e<\e_{0}}\|g(\e,\cd)\|_{C}<\infty$ for a sufficiently small $\e_{0}>0$. Furthermore, for any small $\e>0$
\ali
{
|g(\e,\om)-g(\e,\up)|\leq c_{\adr{bdge2}}e^{\theta c_{\adr{bde.func}}}c_{\adr{bde.func}}d_{\theta}(\om,\up)
}
for $\om,\up\in X$ with $\om_{0}=\up_{0}$. This yields the boundedness $\sup_{0<\e<\e_{0}}\|g(\e,\cd)\|_{F_{\theta}}<\infty$.
\smallskip
\par
Since $\lam$ is in the resolvent of $\RR\in \LR(F_{\theta,b}(X))$, we have $\|\RR_{\lam}\|_{F_{\theta}\to L^{1}(\nu)}\leq \|\RR_{\lam}\|_{F_{\theta}}<\infty$. Now we show that $\RR_{\lam}\,:\,F_{\theta,b}(X)\to L^{1}(\nu)$ is weak bounded. Choose any $f_{\e},f\in F_{\theta,b}(X)$ so that $\sup_{\e}\|f_{\e}\|_{F_{\theta}}<\infty$ and $\|f_{\e}-f\|_{L^{1}}\to 0$. Letting $f_{\e}:=f_{\e}-f$, we may assume $\|f_{\e}\|_{L^{1}(\nu)}\to 0$. 
By $\|\RR_{\lam} f_{\e}\|_{C}+[\RR_{\lam} f_{\e}]_{\theta}^{1}= \|\RR_{\lam} f_{\e}\|_{F_{\theta}}\leq c_{\adl{bRlfe}}:=\sup\|f_{\e}\|_{F_{\theta}}\|\RR_{\lam}\|_{F_{\theta}}<\infty$, we see that  $\{\RR_{\lam}f_{\e}\}_{\e}$ is equicontinuous and uniformly bounded. By virtue of Ascoli Theorem again, there exist a subsequence $(\e_{k})$ and $g\in C_{b}(X,\R)$ such that $\RR_{\lam} f_{\e_{k}}$ converges to $g$ for each point in $X$. Now we prove $g=0$. By Lebesgue dominated convergence theorem, $\SR f_{\e^\p_{k}}$ converges to $g$ in $L^{1}(\nu)$ for some subsequence $(\e^\p_{k})$ of $(\e_{k})$. We have
\ali
{
\|(\RR-\lam\IR)g\|_{L^{1}(\nu)}\leq&\|(\RR-\lam\IR)g-(\RR-\lam\IR)\RR_{\lam} f_{\e^\p_{k}}\|_{L^{1}(\nu)}+\|(\RR-\lam\IR)\RR_{\lam} f_{\e^\p_{k}}\|_{L^{1}(\mu)}\\
\leq &\|\RR-\lam\IR\|_{L^{1}(\nu)}\|g-\RR_{\lam} f_{\e^\p_{k}}\|_{L^{1}(\nu)}+\|f_{\e^\p_{k}}\|_{L^{1}(\nu)}\to 0
}
as $k\to \infty$. Thus $(\RR-\lam\IR)g=0$ $\nu$-a.e. and $(\RR-\lam\IR)g=0$ in $X$ by continuity. Since $\lam$ is in the resolvent of $\RR\in \LR(F_{\theta,b}(X))$, we get $g=0$. This implies $\RR_{\lam} f_{\e}\to 0$ in $L^{1}(\nu)$ and therefore $\RR_{\lam}$ is weak bounded by Proposition \ref{prop:prop_weakbdd}. Finally, the condition (c) of Theorem \ref{th:asymp_efunc_nB_gene} is obtained by (\ref{eq:tLke->0}).
By virtue of this theorem, we obtain the assertions of this claim.
\cla
{\label{cla:tnune_Fthe->0}
$\|\ti{\nu}_{n}(\e,\cd)\|_{F_{\theta}\to \K}\to 0$.
}
We put
\ali
{
\BR^{0}=C_{b}(X),\quad \DR^{0}=\BR^{1}=\cdots=\BR^{n+1}=F_{\theta,b}(X).
}
Then the all conditions (a)-(c) in Theorem \ref{th:asymp_eval_nB_gene3} are guaranteed by a similar argument above. This theorem yields the assertion of this claim together with the fact $|\nu(\e,h)|\leq \|h\|_{C}<+\infty$.
\cla
{
$\|\ti{h}_{n}(\e,\cd)\|_{L^{1}(\nu)}\to 0$ and $\sup_{\e>0}\|\ti{h}_{n}(\e,\cd)\|_{F_{\theta}}<+\infty$.
}
Recalling the form (\ref{eq:thne=...}) of $\ti{h}_{n}(\e,\cd)$, we show $\ti{c}_{n}(\e)\to 0$ and $\nu(h(\e,\cd))\to 1$. To do this, we firstly check $\nu_{k}(\ti{g}_{n-k}(\e,\cd))\to 0$ in the form (\ref{eq:tcne=...}) of $\ti{c}_{n}(\e)$. The boundedness of $\SR$, $\TR(\e,\cd)$, $(\tLR_{k}(\e,\cd))$ and $g(\e,\cd)$ in $F_{\theta,b}(X)$ imply that the function $\ti{g}_{k}(\e,\cd)$ in the form (\ref{eq:tgke=...in2}) is also bounded with respect to the norm $\|\cd\|_{F_{\theta}}$. In addition to $\|\ti{g}_{k}(\e,\cd)\|_{L^{1}(\nu)}\to 0$, this fact yields the convergence $\kappa_{k}(\ti{g}_{n-k}(\e,\cd))\to 0$ by the weak boundedness of $\LR_{i}$ and $\RR_{\lam}$. Therefore, it follows from the form (\ref{eq:nuk=...}) that $\nu_{k}(\ti{g}_{n-k}(\e,\cd))\to 0$. We secondly prove $\ti{\nu}_{n}(\e,g(\e,\cd))\to 0$. This is true by the estimate $|\ti{\nu}_{n}(\e,g(\e,\cd))|\leq \|\ti{\nu}_{n}(\e,\cd)\|_{F_{\theta}\to \K}\|g(\e,\cd)\|_{F_{\theta}}$ and the claim \ref{cla:tnune_Fthe->0}. Thus we see $\ti{c}_{n}(\e)\to 0$. Finally, we will show $\nu(h(\e,\cd))\to 1$. Note the form $\nu(\e,g(\e,\cd))=1/\nu(h(\e,\cd))$. We have
$\nu(\e,g(\e,\cd))=\nu(g(\e,\cd))+\ti{\nu}_{0}(\e,g(\e,\cd))\to \nu(h)=1$ and then $\nu(h(\e,\cd))\to 1$. Consequently, we obtain convergence $\|\ti{h}_{n}(\e,\cd)\|_{L^{1}(\nu)}\to 0$. The boundedness of $\|\ti{h}_{n}(\e,\cd)\|_{F_{\theta}}$ is guaranteed by the boundedness of $\|g_{i}\|_{F_{\theta}}$ and $\|\ti{g}_{n}(\e,\cd)\|_{F_{\theta}}$. 
\cla
{
$\ti{\mu}_{n}(\e,f)\to 0$ for each $f\in F_{\theta,b}(X)$. 
}
Recall the form (\ref{eq:tnumef=...}) of $\ti{\mu}_{n}(\e,f)$. Noting $\|\ti{h}_{n}(\e,\cd)f\|_{L^{1}(\nu)}\to 0$ and $\sup_{\e>0}\|\ti{h}_{n}(\e,\cd)f\|_{F_{\theta}}<+\infty$, convergence $\nu_{k}(\ti{h}_{n-k}(\e,\cd)f)\to 0$ is satisfied by a similar proof for $\nu_{k}(\ti{g}_{n-k}(\e,\cd))\to 0$. Moreover, $\ti{\nu}_{n}(\e,h(\e,\cd)f)\to 0$ by the claim \ref{cla:tnune_Fthe->0} and $\sup_{\e>0}\|h(\e,\cd)f\|_{F_{\theta}}<+\infty$. Hence the assertion is complete.
\proe
\rem
{
Assume that $X$ is a topological Markov shift whose transition matrix is finitely irreducible and whose state space is finite. In this case, if $\sup_{\e>0}\|\ti{\ph}_{n}(\e,\cd)\|_{F_{\theta}}<+\infty$ then $\sup_{\e>0}\|\tLR_{n}(\e,\cd)\|_{F_{\theta}}<+\infty$ \cite{T2011}.
}
\thms
\label{th:asymp_mue2}
Let $X$ be a topological Markov shift whose transition matrix is finitely irreducible. Assume that $\ph,\ph_{1},\dots, \ph_{n}$ satisfy the condition ($\Ph$). Assume also convergence $(\|\SR\|_{F_{\theta(\e)}})^{n-k+1}\|\tLR_{k}(\e,\cd)\|_{F_{\theta(\e)}}\to 0$ as $\e\to 0$ for each $k=0,1,\dots, n$. Then $\|\ti{h}_{n}(\e,\cd)\|_{C}\to 0$ and $\ti{\mu}_{n}(\e,f)\to 0$ for each $f\in F_{\theta,b}(X)$.
\thme
\pros
If $\dim F_{\theta(\e),b}(X)=1$ then the assertion is trial. Therefore we assume $\dim F_{\theta(\e),b}(X)\geq 2$. In this case, we have $\inf_{\e>0}\|\SR\|_{\BR^{i}_{\e}}>0$ (see Claim \ref{cla:pos_S} in the proof of Theorem \ref{th:asymp_efunc_nospecgap}). This means that the assumption also implies $\|\tLR_{n}(\e,\cd)\|_{C}\leq \|\tLR_{n}(\e,\cd)\|_{F_{\theta(\e)}}\to 0$. Thus $\ti{\nu}_{n}(\e,f)\to 0$ for each $f\in F_{\theta,b}(X)$ by Theorem \ref{th:asymp_pphe}.
we set
\ali
{
\DR^{0}_{\e}=\BR^{0}_{\e}=\cdots=\BR^{n+1}_{\e}=F_{\theta(\e),b}(X).
}
Now we check the conditions of Theorem \ref{th:asymp_mue2}. The condition (a) follows from $\LR(F_{\theta,b}(X))\subset \LR(F_{\theta(\theta),b}(X))$. By 
$|\nu(f)|\leq \|f\|_{C}\leq \|f\|_{F_{\theta(\e)}}$ for $f\in F_{\theta(\e),b}(X)$ and $\|h\|_{F_{\theta(\e)}}\leq \|h\|_{F_{\theta}}$, the condition (b) is valid. Finally the condition (c) is obtained from the assumption. As a result, Theorem \ref{th:asymp_mue2} says $\|\ti{g}_{n}(\e,\cd)\|_{C}\to 0$ from $\|\ti{g}_{n}(\e,\cd)\|_{C}\leq \|\ti{g}_{n}(\e,\cd)\|_{F_{\theta(\e)}}$. To see convergence of $\ti{h}_{n}(\e,\cd)$, we recall the form (\ref{eq:tcne=...2}) of $\ti{c}_{n}(\e)$ in the expression (\ref{eq:thne=...}) of $\ti{h}_{n}(\e,\cd)$. Then we have the estimate
\ali
{
|\ti{c}_{n}(\e)|\leq \sum_{k=0}^{n}\ti{\nu}_{n-k}(\e,g_{k})+\|\ti{g}_{n}(\e,\cd)\|_{C}
}
and this vanishes as $\e\to 0$ by $\ti{\nu}_{n-k}(\e,g_{k})\to 0$ and $\|\ti{g}_{n}(\e,\cd)\|_{C}\to 0$. Since $\ti{c}_{n}(\e)$ is the remainder part of $1/\nu(h(\e,\cd))$, we get $\|\ti{h}_{n}(\e,\cd)\|_{C}\to 0$. Finally, we obtain the expansion
\ali
{
\mu(\e,f)=\nu(\e,h(\e,\cd)f)=&\sum_{k=0}^{n}\nu(\e,h_{k}f)\e^{k}+\nu(\e,\ti{h}_{n}(\e,\cd)f)\e^{n}\\
=&\sum_{k=0}^{n}\sum_{i=0}^{n-k}\nu_{i}(h_{k}f)\e^{k+i}+(\sum_{k=0}^{n}\ti{\nu}_{k}(\e,h_{n-k}f)+\nu(\e,\ti{h}_{n}(\e,\cd)f))\e^{n}
}
and the remainder part vanishes by $\ti{\nu}_{k}(\e,h_{n-k}f)\to 0$ and $|\nu(\e,\ti{h}_{n}(\e,\cd)f)|\leq \|\ti{h}_{n}(\e,\cd)f\|_{C}$.
\proe
\subsection{Asymptotic perturbation of Gibbs measure associated with graph directed Markov systems}\label{sec:GDMS_Gibbs}
In this section, we consider asymptotic behaviours of the Gibbs measure associated with the Hausdorff dimensions of the limit sets of perturbed graph directed Markov systems introduced by \cite{MU}.
A set $G=(V,E,i(\cd),t(\cd))$ is called a directed multigraph if this consists of a countable vertex set $V$, a countable edge set $E$, and two maps $i(\cd)$ and $t(\cd)$ from $E$ to $V$. For each $e\in E$, $i(e)$ is called the initial vertex of $e$ and $t(e)$ called the terminal vertex of $e$. Denoted by $X$ the one-sided shift space (see (\ref{eq:X=...})) endowed with the countable state space $E$ and with the incidence matrix $A=(A(ee^\p))_{E\times E}$ which is defined by
$A(ee^\p)=1$ if $t(e)=i(e^\p)$ and $A(ee^\p)=0$ otherwise.
We will use the notation $\si, d_{\theta}, [f]_{\theta}^{n}, \|f\|_{C}, \|f\|_{F_{\theta}}, C_{b}(X), F_{\theta,b}(X), P(\ph)$ defined in Section \ref{sec:thermo}.
We begin with the definition of this system. Let $D$ be a positive integer, $\beta\in (0,1]$ and $r\in (0,1)$. We introduce a set $(G,(J_{v}),(O_{v}),(T_{e}))$ with $\sharp V<\infty$ satisfying the following conditions (i)-(iv):
\ite
{
\item[(i)] For each $v\in V$, $J_{v}$ is a compact and connected subset of $\R^{D}$ satisfying that the interior $\mathrm{int} J_{v}$ of $J_{v}$ is not empty, and $\mathrm{int} J_{v}$ and $\mathrm{int}J_{v^\p}$ are disjoint for $v^\p\in V$ with $v\neq v^\p$.
\item[(ii)] For each $v\in V$, $O_{v}$ is a bounded, open and connected subset of $\R^{D}$ containing $J_{v}$.
\item[(iii)] For each $e\in E$, a function $T_{e}\,:\,O_{t(e)}\to T_{e}(O_{t(e)})\subset O_{i(e)}$ is a $C^{1+\beta}$-conformal diffeomorphism with $T_{e}(\mathrm{int}J_{t(e)})\subset \mathrm{int}J_{i(e)}$ and $\sup_{x\in O_{t(e)}}\|T_{e}^\p(x)\|\leq r$, where $\|T_{e}^\p(x)\|$ means the operator norm of $T_{e}^\p(x)$. Moreover, for any $e,e^\p\in E$ with $e\neq e^\p$ and $i(e^\p)=i(e)$, $T_{e}(\mathrm{int}J_{t(e)})\cap T_{e^\p}(\mathrm{int}J_{t(e^\p)})=\emptyset$, namely the open set condition (OSC) is satisfied.
\item[(vi)] (Bounded distortion) There exists a constant $c_{\adl{Gbd}}>0$ such that for any $e\in E$ and $x,y\in O_{t(e)}$, $|\|T_{e}^\p(x)\|-\|T_{e}^\p(y)\||\leq c_{\adr{Gbd}}\|T_{e}^\p(x)\| |x-y|^{\beta}$, where $|\cd|$ means a norm of any Euclidean space. 
\item[(v)] (Cone condition) If $\sharp E=\infty$ then there exist $\gamma, l>0$ with $\gamma<\pi/2$ such that for any $v\in V$, $x\in J_{v}$, there is $u\in \R^{D}$ with $|u|=1$ so that the set $\{y\in \R^{D}\,:\,0<|y-x|<l \text{ and }(y-x,u)>|y-x|\cos \gamma\}$ is in $\mathrm{int} J_{v}$, where $(y-x,u)$ denotes the inner product of $y-x$ and $u$.
}
Under these conditions (i)-(v), we call the set $(G,(J_{v}), (O_{v}), (T_{e}))$ a {\it graph directed Markov system} ({\it GDMS} for short).
The Hausdorff dimension of the limit set of this system has been mainly studied by many authors \cite{MU1999,MU,RU,RU2}.
\smallskip
\par
The coding map $\pi\,:\,E^{\infty}\to \R^{D}$ is defined by $\pi\om=\bigcap_{n=0}^{\infty}T_{\om_{0}}\cdots T_{\om_{n}}(J_{t(\om_{n})})$ for $\om \in E^{\infty}$. Put $K=\pi(E^{\infty})$. This set is called the {\it limit set} of the GDMS.
We define a function $\ph\,:\,E^{\infty}\to \R$ by
\ali
{
\ph(\om)=(1/D)\log |g(\om)|\quad \text{ with }\quad g(\om)=\det(T_{\om_{0}}^\p(\pi\si\om)).
}
Note that the function $g(\om)$ is used in the proof of Lemma \ref{lem:asymp_Lphe}. Remark also that if the graph $G$ is infinite, then $\|\ph\|_{C}=+\infty$ holds.
Put
\ali
{
\underline{s}=\inf\{s\geq 0\,:\,P(s\ph)<+\infty\},
}
where $P(s\ph)$ means the topological pressure of $s\ph$ which is given by (\ref{eq:P(ph)=}). We call the GDMS {\it regular} if $P(s\ph)=0$ for some $s\geq \underline{s}$. 
The GDMS is said to be {\it strongly regular} if $0<P(s\ph)<+\infty$ for some $s\geq \underline{s}$ (see \cite{MU,RU} for the terminology).
It is known that the general Bowen's formula is satisfied:
\thm
{[\cite{RU}]\label{th:MU}
Let $(G,(J_{v}),(O_{v}),(T_{e}))$ be a graph directed Markov system. Assume that $E^{\infty}$ is finitely irreducible. Then $\dim_{H}K=\inf\{t\in \R\,:\,P(t\ph)\leq 0\}$. In addition to the above condition, we also assume that the potential $\ph$ is regular. Then $s=\dim_{H}K$ if and only if $P(s\ph)=0$.
}
Now we formulate an asymptotic perturbation of graph directed Markov systems. Fix integers $n\geq 0$, $D\geq 1$ and a number $\beta\in (0,1]$. Consider the following conditions $(G.1)_{n}$ and $(G.2)_{n}$:
\ite
{
\item[$(G.1)_{n}$] The code space $E^{\infty}$ is finitely irreducible. The set $(G,(J_{v}),(O_{v}),(T_{e}))$ is a GDMS on $\R^{D}$ with strongly regular and the limit set $K$ has positive dimension. Moreover, the function $T_{e}$ is of class $C^{1+n+\beta}(O_{t(e)})$ for each $e\in E$.
\item[$(G.2)_{n}$] The set $\{(G,(J_{v}),(O_{v}),(T_{e}(\e,\cd)))\,:\,\e>0\}$ is a GDMS with a small parameter $\e>0$ satisfying the following (i)-(vi):
\ite
{
\item[(i)] For each $e\in E$, the function $T_{e}(\e,\cd)$ has the $n$-asymptotic expansion:
\ali
{
T_{e}(\e,\cd)=T_{e}+T_{e,1}\e+\cdots+T_{e,n}\e^{n}+\ti{T}_{e,n}(\e,\cd)\e^{n}\ \text{ on }J_{t(e)}
}
for some functions $T_{e,k}\in C^{1+n-k+\beta}(O_{t(e)},\R^{D})$ $(k=1,2,\dots, n)$ and $\ti{T}_{e,n}(\e,\cd)\in C^{1+\beta(\e)}(O_{t(e)},\R^{D})$ $(\beta(\e)>0)$ satisfying $\sup_{e\in E}\sup_{x\in J_{t(e)}}|\ti{T}_{e,n}(\e,x)|\to 0$.
\item[(ii)] There exist constants $t(l,k)\in (0,1]$ ($l=0,1,\dots, n$,\ $k=1,\dots, n-l+1$) such that the function $x\mapsto T_{e,l}^{(k)}(x)/\|T_{e}^{\p}(x)\|^{t(l,k)}$ is bounded, $\beta$-H\"older continuous and its H\"older constant is bounded uniformly in $e\in E$.
\item[(iii)] $c_{\adl{Gbd3}}(\e):=\sup_{e\in E}\sup_{x\in J_{t(e)}}(\|\frac{\partial}{\partial x}\ti{T}_{e,n}(\e,x)\|/\|T_{e}^\p(x)\|^{\ti{t}_{0}})\to 0$ as $\e\to 0$ for some $\ti{t}_{0}\in (0,1]$.
\item[(iv)] $\dim_{H}K/D>p(n)$, where $p(n)$ is taken by
\alil
{
p(n):=&
\case
{
\underline{p}/\ti{t},&n=0\\
\max\big(\underline{p}+n(1-t_{1}), \underline{p}+n(1-t_{2})/2, \cdots, \underline{p}+n(1-t_{n})/n,\\
\qqqqqqquad\underline{p}/t_{1},\ \underline{p}/t_{2},\ \cdots, \underline{p}/t_{n},\ \underline{p}+1-\ti{t},\ \underline{p}/\ti{t}\big),& n\geq 1.
}\\
t_{k}:=&\min\{\frac{1}{D}\sum_{p=1}^{D}t(i_{p},j_{p}+1)\,:\,i:=i_{1}+\cdots+i_{D}\text{ and }j:=j_{1}+\cdots+j_{D}\text{ satisfy}\label{eq:tk=}\\
&\qqqqqquad i=k \text{ and }j=0 \text{ or }0\leq i<k \text{ and }1\leq i+j\leq k\}\nonumber\\
\ti{t}:=&\min\left\{t_{n},\ \ti{t}_{0},\ \frac{\ti{t}_{0}}{D}+\frac{D-1}{D}t(1,1),\dots,\ \frac{\ti{t}_{0}}{D}+\frac{D-1}{D}t(n,1)\right\}\label{eq:tt=}\\
\underline{p}:=&\underline{s}/D.\nonumber
}
}
}
Note that if the edge set $E$ is finite, then the conditions (ii) and (iv) are always satisfied because $\|T_{e}^\p(x)\|$ is uniformly bounded away from zero, and $p(n)$ becomes zero by taking $t(l,k)\equiv 1$. Moreover, $c_{\adr{Gbd3}}(\e)$ in (iii) can be taken as $\sup_{e\in E}\sup_{x\in J_{t(e)}}\|\frac{\partial }{\partial x}\ti{T}_{e,n}(\e,x)\|$ when $E$ is finite. Let $K(\e)$ be the limit set of the perturbed GDMS $(G,(J_{v}),(O_{v}),(T_{e}(\e,\cd)))$. We put $s(\e)=\dim_{H}K(\e)$. Under those conditions, we obtained the following:
\thm
{[\cite{T2020pre}]\label{th:asymp_dim}
Assume that the conditions $(G.1)_{n}$ and $(G.2)_{n}$ are satisfied with fixed integer $n\geq 0$. Then the perturbed GDMS $(G,(J_{v}),(O_{v}),(T_{e}(\e,\cd)))$ is strongly regular for any small $\e>0$, and 
there exist $s_{1},\dots, s_{n}\in \R$ such that the Hausdorff dimension $s(\e)$ of the limit set $K(\e)$ of the perturbed system has the form $s(\e)=s(0)+s_{1}\e+\cdots+s_{n}\e^{n}+o(\e^{n})$ as $\e \to 0$ with $s(0)=\dim_{H}K$.
}
\rem
{
Roy and Urba\'nski \cite{RU} considered continuous perturbation of infinitely conformal iterated function systems given as a special GDMS. They also studied analytic perturbation of GDMS with $D\geq 3$ in \cite{RU2}. We investigated an asymptotic perturbation of GDMS with finite graph in \cite{T2016}. Theorem \ref{th:asymp_dim} is an infinite graph version of this previous result in \cite{T2016}
}
In order to investigate an asymptotic perturbation of the Gibbs measure associated with the dimension $\dim_{H}K(\e)$, 
we further introduce the following condition:
\ite
{
\item[$(G.3)_{n}$] 
$\di\sup_{\e>0}\sup_{\e\in E}\sup_{x,y\in O_{t(e)}\,:\,x\neq y}\frac{|\frac{\partial}{\partial  x}\ti{T}_{e,n}(\e,x)-\frac{\partial}{\partial  x}\ti{T}_{e,n}(\e,y)|}{|x-y|^{\beta}}<\infty$.
}
Note that such a condition is firstly given in \cite[the condition $(G)^\p_{n}$]{T2016}. We denote the physical potential for the perturbed GDMS by
\alil
{
\ph(\e,\om)=(1/D)\log |g(\e,\om)|\quad \text{ with }\quad g(\e,\om)=\textstyle{\det(\frac{\partial}{\partial x}T_{\om_{0}}(\e,\pi(\e,\si\om)))}\label{eq:phe=}
}
for $\om\in E^{\infty}$, where $\pi(\e,\cd)$ means the coding map of $K(\e)$ which is defined by $\pi(\e,\om)=\bigcap_{n=0}^{\infty}T_{\om_{0}}(\e,\cd)\cdots T_{\om_{n}}(\e,J_{t(\om_{n})})$ for $\om\in E^{\infty}$. The condition $(G.3)_{n}$ implies the boundedness of the Lipschitz constant of the remainder $\ti{\ph}_{n}(\e,\cd)$ of $\ph(\e,\cd)$.
\smallskip
\par
Now we are in a position to state the result in this section.
\thm
{\label{th:asymp_Gibbs_GDMS}
Assume that the conditions $(G.1)_{n}$-$(G.3)_{n}$ are satisfied with fixed integer $n\geq 0$. Then there exist bounded functionals $\mu_{1},\dots, \mu_{n}$ in the dual $(F_{\theta,b}(E^{\infty},\C))^{*}$ of $F_{\theta,b}(E^{\infty},\C)$
 
 such that the Gibbs measure $\mu(\e,\cd)$ of the potential $\ph(\e,\cd)$ satisfies the asymptotic expansion
\alil
{
\mu(\e,f)=&\mu(f)+\mu_{1}(f)\e+\cdots+\mu_{n}(f)\e^{n}+o(\e^{n})\label{eq:muef=} 
}
as $\e\to 0$ for $f\in F_{\theta,b}(E^{\infty},\C)$, where $\mu(\e,\cd)$ is the Gibbs measure of $(\dim_{H}K(\e))\phe$ and $\mu$ is the Gibbs measure of $(\dim_{H}K(0))\ph$. 
}
\rem
{
This theorem is an infinite graph version of \cite[Theorem 1.3]{T2016} which considered an asymptotic expansions of the Gibbs measure
 in the case $\sharp E<\infty$.
}
To show the theorem \ref{th:asymp_Gibbs_GDMS}, we start with the following previous lemma:
\lem
{[\cite{T2016,T2020pre}]\label{lem:asymp_pi}
Assume that the conditions $(G.1)_{n}$ and $(G.2)_{n}$ are satisfied. Choose any $r_{1}\in (r,1)$. Then there exist functions $\pi_{1},\pi_{2},\dots, \pi_{n}\in F_{r_{1},b}(E^{\infty},\R^{D})$ and $\ti{\pi}(\e,\cd)\in C_{b}(E^{\infty},\R^{D})$ such that
$\pi(\e,\cd)=\pi+\pi_{1}\e+\cdots+\pi_{n}\e^{n}+\ti{\pi}_{n}(\e,\cd)\e^{n}$ and $\|\ti{\pi}_{n}(\e,\cd)\|_{C}:=\sup_{\om\in E^{\infty}}|\ti{\pi}_{n}(\e,\om)|\to 0$ as $\e\to 0$.
}
\pros
See \cite[Lemma 3.13]{T2020pre} and \cite[Lemma 3.1]{T2016}.
\proe
\lem
{\label{lem:bd_asymp_pi}
Assume that the conditions $(G.1)_{n}$-$(G.3)_{n}$ are satisfied. Then we have $\limsup_{\e\to 0}[\ti{\pi}_{n}(\e,\cd)]_{r_{2}}^{1}<\infty$ with any number $r_{2}\in (r_{1},1)$.
}
\pros
This follows from \cite[Lemma 3.3]{T2016} replacing $\max_{e\in E}$ by $\sup_{e\in }E$. Remark that $\max_{v\in V}$ in this lemma need not be replaced with $\sup_{v\in V}$ by the assumption $\sharp V<\infty$.
\proe
We define an operator $\LR_{k}\,:\,F_{\theta,b}(X)\to F_{\theta,b}(X)$ by
\ali
{
\LR_{k}f=\LR_{s(0)\ph}(G_{k} f)
}
with
\ali
{
G_{k}=\di\sum_{q=0}^{k}\frac{s_{q,k}}{q!}\ph^{q}+\sum_{v=1}^{k}\sum_{q=0}^{k-v}\sum_{l=0}^{v}\sum_{j=0}^{\min(l,q)}\frac{s_{q,k-v}\cd a_{l,j}}{(q-j)!}\ph^{q-j}l!\sum_{j_{1},\dots,j_{v}\geq 0\,:\atop{j_{1}+\cdots+j_{v}=l\atop{j_{1}+2j_{2}+\dots+v j_{v}=v}}}\prod_{u=1}^{v}\Big(\frac{(g_{u})^{j_{u}}}{j_{u}!g^{j_{u}}}\Big)
}
where $\LR_{0}f$ means $\sum_{e\in E\,:\,t(e)=i(\om_{0})}f(e\cd\om)$ and where $a_{l,j}$ and $s_{q,k-v}$ are given by expanding $\binom{p}{l}=a_{l,0}+\sum_{j=1}^{l}a_{l,j}(p-s)^{j}$ with $a_{l,0}=\binom{s}{l}$ of binomial coefficient and $(s(\e)-s(0))^{k}=\sum_{i=0}^{n}s_{k,i}\e^{i}+o(\e^{n})$, respectively.
Here $g,g_{1},\dots, g_{n},g(\e,\cd)$ satisfy the equation $g(\e,\cd)=g+\sum_{k=1}^{n}g_{k}\e^{k}+\ti{g}_{n}(\e,\cd)\e^{n}$. 
We have the following:
\lem
{\label{lem:asymp_Lphe}
Assume that the conditions $(G.1)_{n}$-$(G.3)_{n}$ are satisfied. Then each $\LR_{k}$ is a bounded linear operator acting on $F_{\theta,b}(X)$ with $\theta=r_{2}^{\beta}$ and the Ruelle operator $\LR(\e,\cd)$ of the potential $s(\e)\phe$ has the expansion
$\LR(\e,\cd)=\LR_{s\ph}+\sum_{k=1}^{n}\LR_{k}\e^{k}+\tLR_{n}(\e,\cd)\e^{n}$
with $\|\ti{\LR}_{n}(\e,\cd)\|_{C}\to 0$. 
}
\pros
This equation is guaranteed by the form (3.37) in \cite{T2020pre}, and convergence of $\ti{\LR}_{n}(\e,\cd)$ is yielded by the proof of \cite[Theorem 1.1]{T2020pre}.
\proe
\lem
{\label{lem:bd_asymp_phe}
Assume that the conditions $(G.1)_{n}$-$(G.3)_{n}$ are satisfied. Then we have $\sup_{\e>0}[\ti{\ph}(\e,\cd)]_{\theta}^{1}<\infty$ with $\theta=r_{2}^{\beta}$, where $\phe$ appears in Lemma \ref{lem:asymp_Lphe} and $r_{2}$ is given in Lemma \ref{lem:bd_asymp_pi}. In particular, $\sup_{\e>0}\|\ti{\LR}_{n}(\e,\cd)\|_{F_{\theta}}<\infty$.
}
\pros
Firstly we show the former assertion $\sup_{\e>0}[\ti{\ph}_{n}(\e,\cd)]_{\theta}^{1}<\infty$. Denote by $\ti{T}_{e,n}(\e,\cd)=(\ti{t}_{e,n,1}(\e,\cd),\cdots, \ti{t}_{e,n,D}(\e,\cd))$. By virtue of the condition $(G.3)_{n}$, we see
\ali
{
\sup_{\e>0}\sup_{\e\in E}\sup_{x,y\in O_{t(e)}\,:\,x\neq y}\frac{|\frac{\partial}{\partial  x_{j}}\ti{t}_{e,n,i}(\e,x)-\frac{\partial}{\partial  x_{j}}\ti{t}_{e,n,i}(\e,y)|}{|x-y|^{\beta}}<\infty
}
for each $i,j=1,\dots, D$. This yields the boundedness of the $\beta$-H\"older constant of the remainder part of $\det T^\p_{e}(\e,\cd)$ uniformly in $e\in E$. Thus the $\theta$-Lipschitz constant of the remainder part of $\det T^\p_{\om_{0}}(\e,\pi(\e,\si\om))$ is also finite. By using the form \cite[(3.14)]{T2016} of $\ti{\ph}_{n}(\e,\cd)$, we obtain the boundedness of $[\ti{\ph}_{n}(\e,\cd)]_{\theta}^{1}$ uniformly in $\e>0$ (see also \cite[Lemma 3.4]{T2016}).
Secondly we check $\sup_{\e>0}\|\ti{\LR}_{n}(\e,\cd)\|_{F_{\theta}}<\infty$.
To do this, we prove that $\ti{\LR}_{n}(\e,\cd)$ is defined as (\ref{eq:tLRke=}) replacing $\phe:=s(\e)\phe$, $\ph:=s(0)\ph$ and $\ph_{k}:=\sum_{i=0}^{k}s_{i}\ph_{k-i}$. Note that this $\ph_{k}$ is the $k$-th coefficient of the expansion of $s(\e)\phe$. This is shown by checking the equation $G_{k}=F_{k}$, where $F_{k}$ is given in (\ref{eq:Fk=}) replacing by $\ph_{k}:=\sum_{i=0}^{k}s_{i}\ph_{k-i}$. Fix $\om\in X$. By Lemma \ref{lem:asymp_Lphe}, the function $e^{s(\e)\ph(\e,\om)}=\LR(\e,\chi_{[\om_{0}]})(\si\om)$ has the expansion
\alil
{
e^{s(\e)\ph(\e,\om)}=e^{s(0)\ph(\om)}+\sum_{k=1}^{n}e^{s(0)\ph(\om)}G_{k}(\om)\e^{k}+\ti{G}_{n}(\e,\om)\e^{n}\label{eq:|ge|^se=...}
}
with $\ti{G}_{n}(\e,\om):=\tLR_{n}(\e,\chi_{[\om_{0}]})(\si\om)\to 0$ as $\e\to 0$, where $\chi_{[\om_{0}]}$ is the indicator on the cylinder set $[\om_{0}]=\{\up\in X\,:\,\up_{0}=\om_{0}\}$. Therefore, $s(\e)\ph(\e,\om)=\log(e^{s(\e)\ph(\e,\om)})$ is also asymptotically expanded. On the other hand, (\ref{eq:Fk=}) and (\ref{eq:tFke=}) imply
\alil
{
e^{s(\e)\ph(\e,\om)}=e^{s(0)\ph(\om)}+\sum_{k=1}^{n}\e^{s(0)\ph(\om)}F_{k}(\om)\e^{k}+\e^{s(0)\ph(\om)}\ti{F}_{n}(\e,\om)\e^{n}\label{eq:e^sephe=...}
}
with $\ti{F}_{n}(\e,\om)\to 0$ as $\e\to 0$. Thus the equations (\ref{eq:|ge|^se=...}) and (\ref{eq:e^sephe=...}) imply $G_{k}(\om)=F_{k}(\om)$ for $k=1,2,\dots, n$ by induction. Consequently, $\tLR_{n}(\e,\cd)$ is defined by (\ref{eq:tLRke=}). The finiteness $\sup_{\e>0}[\ti{\ph}_{n}(\e,\cd)]_{\theta}^{1}<\infty$ and $[\ph_{k}]_{\theta}^{1}<\infty$ yield $\sup_{\e>0}[\ti{\ph}_{i}(\e,\cd)]_{\theta}^{1}<\infty$ for any $0\leq i\leq n$ and this is able to show $[\ti{F}_{n}(\e,\om)]_{\theta}^{1}<\infty$ by the form (\ref{eq:tFke=}). Hence the $\theta$-Lipschitz norm of $\tLR_{n}(\e,\cd)=\LR_{s(0)\ph}(\ti{F}_{n}(\e,\cd)\cd)$ is bounded uniformly in $\e>0$.
\proe
\noindent
({\it Proof of Theorem \ref{th:asymp_Gibbs_GDMS}}). The convergence of the reminder $\ti{\mu}_{n}(\e,f)$ of the perturbed Gibbs measure $\mu(\e,f)$ of $f\in F_{\theta,b}(E^{\infty},\R)$ is guaranteed by Theorem \ref{th:asymp_pphe_he} together with Lemma \ref{lem:asymp_Lphe} and Lemma \ref{lem:bd_asymp_phe}. 
\qed
\subsection{Example of perturbation without uniform spectral gap}\label{sec:ex_withoutSG}
Let $X$ be a subshift of finite type whose shift is topologically mixing (i.e. topological Markov shift with finite state space and with finitely primitive transitive matrix). Let $\ph\in F_{\theta}(X,\R)$. Recall the reduced resolvent $\SR=(\RR-\lam\IR)^{-1}(\IR-\PR)$ defined in (\ref{eq:ge=,S=,Te=}). For numbers $\theta(\e)\in (0,1)$ with $\theta(\e)\to 1$ as $\e\to 0$, we set a function $c_{\adl{bS}}(\e)$ by
\ali
{
c_{\adr{bS}}(\e):=\frac{e^{26[\ph]_{\theta}^{1}/(1-\theta(\e))}}{(1-\theta(\e))^{4}}.
}
Then $c_{\adr{bS}}(\e)$ tends to $+\infty$ as $\e\to 0$. We will show in Lemma \ref{lem:bd_S} that $\|\SR\|_{F_{\theta(\e)}}=O(c(\e))$ as $\e\to 0$.
Consider a perturbed function $\ph(\e,\cd)\in F_{\theta(\e)}(X,\R)$ having the form
\ali
{
\ph(\e,\cd)=\ph+\ph_{1}\e+\cdots+\ph_{n}\e^{n}+\ti{\ph}_{n}(\e,\cd)\e^{n}
}
for some integer $n\geq 0$ and functions $\ph_{1},\dots,\ph_{n}\in F_{\theta}(X,\R)$. In this case, if we spectrally decompose the Ruelle operator of $\phe$ acting on $F_{\theta(\e)}(X)$ into $\LR(\e,\cd)=\lam(\e)\PR_{\e}+\RR_{\e}$ with the eigenprojection $\PR_{\e}$ of the Perron eigenvalue $\lam(\e)$ and the remainder operator $\RR_{\e}=\LR(\e,\cd)-\lam(\e)\PR_{\e}$, then the distance of the spectrum of $\RR_{\e}$ and the eigenvalue $\lam(\e)$ is bounded by $(1-\theta(\e))\lam(\e)$. Indeed, the spectral radius of  $\RR_{\e}\,:\,F_{\theta(\e)}(X)\to F_{\theta(\e)}(X)$ is not smaller than $\theta(\e)\lam(\e)$ (e.g. \cite[Theorem 1.5(7)]{Baladi}). Therefore this spectral gap vanishes as $\e\to 0$. 
We assume the following:
\ite
{
\item[(B.1)] if $n\geq 1$ then $\theta(\e)\to 1$ satisfies
\ali
{
c_{\adr{bS}}(\e)=O(\e^{-1/(n+2)})
}
as $\e\to 0$. Note that when $n=0$, we arbitrarily choose $\theta(\e)\in (0,1)$ so that $\theta(\e)\to 1$.
\item[(B.2)]
\ali
{
\|\ti{\ph}_{n}(\e,\cd)\|_{F_{\theta(\e)}}=
\case
{
o(c_{\adr{bS}}(\e))& (n=0)\\
o(\e^{1/(n+2)})& (n\geq 1).\\
}
}
}
Under these conditions, the Perron eigenvalue, the Perron eigenfunction and the Perron eigenvector for the Ruelle operator $\LR(\e,\cd)$ have asymptotically converge with order $n$ as $\e\to 0$. Thus we obtain the following:
\thm
{\label{th:ex_conv_Gibbs_nonspec}
Assume that the conditions (B.1) and (B.2) are satisfied. Then we have the $n$-order asymptotic expansions of the topological pressure $P(\phe)$ and of the Gibbs measure $\mu(\e,\cd)$. 
}
In order to prove this proposition, we start with the estimate of the Lipschitz norm of the reduced resolvent:
\lem
{\label{lem:bd_S}
Let $X$ be a subshift of finite type whose shift is topologically mixing. Then for any $\ph\in F_{\theta}(X,\R)$, there exist constants $c_{\adl{bdS}}>0,\e_{0}>0$ such that for any $0<\e<\e_{0}$, $\|\SR\|_{F_{\theta(\e)}}\leq c_{\adr{bdS}}c_{\adr{bS}}(\e)$.
}
\pros
Remark the inequality $\|g\|_{F_{\theta(\e)}}\leq \|g\|_{F_{\theta}}$ for any $g\in F_{\theta}(X)$. 
Take the spectral triplet $(\lam,h,\nu)$ in Ruelle-Perron-Frobenius Theorem (Theorem \ref{th:Spec_gap}). By using the result of \cite[p141]{Stoyanov}, we have that for any $f\in F_{\theta(\e),b}(X)$ and $n\geq 0$
\ali
{
\|\lam^{-n}\LR^{n}f-\nu(f)h\|_{F_{\theta(\e)}}\leq c_{\adr{Sc}}(\e)c_{\adr{Sc2}}(\e)^{n}\|f\|_{F_{\theta(\e)}}
}
with the constants $c_{\adl{Sc}}(\e)>0$ and $c_{\adl{Sc2}}(\e)\in (0,1)$ which are defined by
\ali
{
c_{\adr{Sc}}(\e)=&\frac{1200}{\theta(\e)(1-\theta(\e))^{3}}([\ph]_{\theta}^{1})^{2}(\sharp S)^{9M}e^{18M\|\ph\|_{C}}e^{18[\ph]_{\theta}^{1}\theta(\e)/(1-\theta(\e))}\\
c_{\adr{Sc2}}(\e)=&\Big(1-\frac{1-\theta(\e)}{4e^{8[\ph]_{\theta}^{1}\theta(\e)/(1-\theta(\e))}(\sharp S)^{2M}e^{4M\|\ph\|_{C}}}\Big)^{1/(2M)},
}
where $S$ is the state space, $M=\min\{m\geq 1\,:\,A^{m}>0\}$ and $A$ is the structure matrix of $X$. We also notice 
\ali
{
\SR f=&(\RR-\lam\IR)^{-1}(\IR-\PR)f\\
=&\sum_{n=0}^{\infty}\frac{\RR^{n}(f-\nu(f)h)}{\lam^{n+1}}\\
=&\frac{1}{\lam}(f-\nu(f)h)+\frac{1}{\lam}\sum_{n=1}^{\infty}(\lam^{-n}\LR^{n}f-\nu(f)h).
}
Thus we obtain
\ali
{
\|\SR f\|_{F_{\theta(\e)}}\leq& \frac{1}{\lam}(\|f\|_{F_{\theta(\e)}}+\|h\|_{F_{\theta(\e)}}\|f\|_{C})+\frac{1}{\lam}\sum_{n=1}^{\infty}c_{\adr{Sc}}(\e)c_{\adr{Sc2}}(\e)^{n}\|f\|_{F_{\theta(\e)}}\\
\leq&\frac{1+\|h\|_{F_{\theta}}+2c_{\adr{Sc3}}}{\lam}c_{\adr{bS}}(\e)\|f\|_{F_{\theta(\e)}}
}
for any small $\e>0$ with the constant $c_{\adl{Sc3}}=9600([\ph]_{\theta}^{1})^{2}{(\sharp S)}^{11 M}{
e^{22 M\|\ph\|_{C}-26[\ph]_{\theta}^{1}}}M$, where this constant equals the limit $\lim_{\e\to 0}(c_{\adr{Sc}}(\e)/(1-c_{\adr{Sc2}}(\e)))/c_{\adr{bS}}(\e)$. Hence the assertion is fulfilled by putting $c_{\adr{bdS}}=(1+\|h\|_{F_{\theta}}+2c_{\adr{Sc3}})/\lam$.
\proe
\lem
{\label{lem:bd_tLR}
Let $X$ be a subshift of finite type. Let $\ph(\e,\cd)\in F_{\theta(\e)}(X,\R)$ be functions with the form
$\ph(\e,\cd)=\ph+\sum_{k=1}^{n}\ph_{k}\e^{k}+\ti{\ph}_{n}(\e,\cd)\e^{n}$
for some integer $n\geq 0$ and functions $\ph,\ph_{1},\dots,\ph_{n}\in F_{\theta}(X,\R)$. Then there exists a constant $c_{\adl{bdtLR}}>0$ such that for any $\theta(\e)\in [\theta,1)$ and for each $0\leq k\leq n$
\ali
{
\|\tLR_{k}(\e,\cd)\|_{F_{\theta(\e)}}\leq
\case
{
c_{\adr{bdtLR}}\|\ti{\ph}_{0}(\e,\cd)\|_{F_{\theta(\e)}}& (n=k=0)\\
c_{\adr{bdtLR}}(\e+\|\ti{\ph}_{n}(\e,\cd)\|_{F_{\theta(\e)}}\e^{n-k})& (n\geq 1).\\
}
}
}
\pros
By the form $\ti{\ph}_{k}(\e,\cd)=\ph_{k+1}\e+\cdots+\ph_{n}\e^{n-k}+\ti{\ph}_{n}(\e,\cd)\e^{n-k}$, it is sufficient to show the assertion that $\|\tLR_{k}(\e,\cd)\|_{F_{\theta(\e)}}=O(\|\ti{\ph}_{k}(\e,\cd)\|_{F_{\theta(\e)}})$ as $\e\to 0$.
Recall the form $\tLR_{k}(\e,\cd)$ in (\ref{eq:LRk=}) with $n:=k$. By virtue of the property that $X$ is subshift of finite type, the operator $\LR$ becomes finite summation. Therefore we have the estimate
\alil
{
\|\tLR_{k}(\e,\cd)\|_{C}\leq \|\LR 1\|_{C}\|\ti{F}_{k}(\e,\cd)\|_{C},\qquad \|\ti{F}_{k}(\e,\cd)\|_{C}\leq c_{\adr{btLke}}\|\ti{\ph}_{k}(\e,\cd)\|_{C}\label{eq:tLR_ke<=...}
}
for some constant $c_{\adl{btLke}}>0$ together with the form $\ti{F}_{k}(\e,\cd)$ in (\ref{eq:tFke=}). On the other hand, the basic estimate
\ali
{
[\LR(F f)]_{\theta(\e)}^{1}\leq \|\LR 1\|_{C}(e^{[\ph]_{\theta(\e)}^{1}\theta(\e)^{2}}[\ph]_{\theta(\e)}^{1}\theta(\e)\|F\|_{C}\|f\|_{C}+\theta(\e)[F]_{\theta(\e)}^{1}\|f\|_{C}+\|F\|_{C}[f]_{\theta(\e)}^{1}\theta(\e))
}
for $F,f\in F_{\theta(\e)}(X)$ implies
\alil
{
[\tLR_{k}(\e,f)]_{\theta(\e)}^{1}\leq \|\LR 1\|_{C}(e^{[\ph]_{\theta}^{1}}[\ph]_{\theta}^{1}\|\ti{F}_{k}(\e,\cd)\|_{C}\|f\|_{C}+[\ti{F}_{k}(\e,\cd)]_{\theta(\e)}^{1}\|f\|_{C}+\|\ti{F}_{k}(\e,\cd)\|_{C}[f]_{\theta(\e)}^{1})\label{eq:[tLRke]<=...}
}
using the fact $[\ph]_{\theta(\e)}^{1}\leq [\ph]_{\theta}^{1}$. Hence $\|\tLR_{k}(\e,\cd)\|_{F_{\theta(\e)}}=O(\|\ti{\ph}_{k}(\e,\cd)\|_{F_{\theta(\e)}})$ follows from the inequalities (\ref{eq:tLR_ke<=...}), (\ref{eq:[tLRke]<=...}) and the fact $[\ti{F}_{k}(\e,\cd)]_{\theta(\e)}^{1}=O(\|\ti{\ph}_{k}(\e,\cd)\|_{F_{\theta(\e)}})$. 
\proe
\noindent
({\it Proof of Theorem \ref{th:ex_conv_Gibbs_nonspec}}). We will check the condition of Theorem \ref{th:asymp_pphe}. Since $\|\ti{\ph}_{n}(\e,\cd)\|_{C}$ vanishes as $\e\to 0$ by $\|\ti{\ph}_{n}(\e,\cd)\|_{C}\leq \|\ti{\ph}_{n}(\e,\cd)\|_{F_{\theta(\e)}}\to 0$, we see $\|\tLR_{n}(\e,\cd)\|_{C}\to 0$. Therefore, $P(\ph(\e,\cd))$ has $n$-order asymptotic expansion. On the other hand, for each $0\leq k\leq n-1$, Lemma \ref{lem:bd_tLR} tells us that $\|\tLR_{k}(\e,\cd)\|_{F_{\theta(\e)}}$ is at least equal to $O(\e)$ for any small $\e>0$. In addition to the fact $\|\SR\|_{F_{\theta(\e)}}=O(c_{\adr{bS}}(\e))$ in Lemma \ref{lem:bd_S} and the assumption (B.1), we obtain
\ali
{
\|\SR\|_{F_{\theta(\e)}}^{n-k+1}\|\tLR_{k}(\e,\cd)\|_{F_{\theta(\e)}}=&O(\e^{-(n-k+1)/(n+2)})O(\e)\\
=&O(\e^{(k+1)/(n+2)})=o(1)
}
as $\e\to 0$. Furthermore,
\ali
{
\|\SR\|_{F_{\theta(\e)}}\|\tLR_{n}(\e,\cd)\|_{F_{\theta(\e)}}=&O(\e^{-1/(n+2)})\e+O(\e^{-1/(n+2)})o(\e^{1/(n+2)})\\
=&O(\e^{(n+1)/(n+2)})+o(1)=o(1).
}
Consequently, we get $\max_{0\leq k\leq n}(\|\SR\|_{F_{\theta(\e)}})^{n-k+1}\|\tLR_{k}(\e,\cd)\|_{F_{\theta(\e)}}\to 0$. By virtue of Theorem \ref{th:asymp_mue2}, we see that the Gibbs measure $\mu(\e,f)$ for $\phe$ is asymptotically expanded for each $f\in F_{\theta,b}(X)$. Hence the proof is complete.
\qed
\subsection{Gou\"ezel and Liverani's abstract perturbation theory}\label{sec:GL}
In this section, we state the relationship between the abstract perturbation theory introduced by Gou\"ezel and Liverani \cite{GL} and our main results. We recall the conditions given in \cite[Section 8]{GL}:
\ite
{
\item[(GL.1)] Banach spaces $\BR^{0},\BR^{1},\dots, \BR^{n+1}$ satisfies $\BR^{0}\supset \BR^{1}\supset \cdots \supset \BR^{n+1}$ as linear subspaces.
\item[(GL.2)] An operator $\LR(\e,\cd)$ with a small parameter $\e\in I:=[-\eta,\eta]$ satisfies that there exists constants $c_{\adl{GLc1}}, c_{\adl{GLc2}}>0$ such that for any $\e\in I$, $f\in \BR^{0}$ and $n\geq 1$, $\|\LR(\e,\cd)^{n}f\|_{\BR^{0}}\leq c_{\adr{GLc1}} c_{\adr{GLc2}}^{n}\|f\|_{\BR^{0}}$.
\item[(GL.3)] There exist $c_{\adl{GLc3}}\in (0,c_{\adr{GLc2}})$ and $c_{\adl{GLc4}}>0$ such that for any $\e\in I$, $f\in \BR^{1}$ and $n\geq 1$, $\|\LR(\e,\cd)^{n}f\|_{\BR^{1}}\leq c_{\adr{GLc1}}c_{\adr{GLc3}}^{n}\|f\|_{\BR^{1}}+c_{\adr{GLc1}}c_{\adr{GLc2}}^{n}\|f\|_{\BR^{0}}$.
\item[(GL.4)] For each $j=1,\dots, n$ and $i=j,j+1,\dots, n+1$, $\LR_{j}$ is a bounded linear operator from $\BR^{i}$ to $\BR^{i-j}$.
\item[(GL.5)] Letting $\Delta_{j}(\e,\cd):=\LR(\e,\cd)-\sum_{i=0}^{j-1}\LR_{i}\e^{i}$, there exists a constant $c_{\adl{GLc5}}>0$ such that for any $\e\in I$, $0\leq j\leq n+1$ and $i=j,j+1,\dots, n+1$, $\|\Delta_{j}(\e,\cd)\|_{\BR^{i}\to \BR^{i-j}}\leq c_{\adr{GLc5}}\e^{j}$.
}
For $\rho\in (c_{\adr{GLc3}},c_{\adr{GLc2}})$ and $\delta>0$, put $V_{\delta,\rho}=\{z\in \C\,:\,|z|\geq \rho \text{ and }\dist(z,\sigma(\LR_{0}\,:\,\BR^{i}\to \BR^{i}))\geq \delta\text{ for all }1\leq i\leq s\}$, where $\sigma(\LR_{0})$ denotes the spectrum of $\LR_{0}$. Under these conditions (GL.1)-(GL.5), the asymptotic expansion of the resolvent $(z\IR-\LR(\e,\cd))^{-1}\,:\,\BR^{n+1}\to \BR^{0}$ for $z\in V_{\delta,\rho}$ was explicit given \cite[Theorem 8.1]{GL} replacing by $n=s-1$, $\e=t$ and $\LR(\e,\cd)=\LR_{t}$. In particular, the perturbed eigenvalue $\lam(\e)$ of an isolated simple eigenvalue of $\LR_{0}$ and the corresponding eigenprojection $\PR(\e,\cd)$ (eigenvector) are of class $C^{n}$. Therefore, roughly speaking,  the main difference between the perturbation result of \cite{GL} and our main result is whether to get a higher-order derivative or get an asymptotic expansion at $\e=0$.
\smallskip
\par
Now we will show that if the conditions (GL.1)-(GL.5) and the two conditions (A.1)(A.2) below are satisfied, then we have the conditions (I)-(III), the conditions in Theorem \ref{th:asymp_eval_nB_gene2} and the conditions in Theorem \ref{th:asymp_eval_nB_gene3}. We give an additional condition for the unperturbed operator $\LR_{0}$:
\ite
{
\item[(A.1)] There is a nonzero isolated simple eigenvalue $\lam\in \C$ of $\LR_{0}\,:\,\BR^{i}\to \BR^{i}$ for all $1\leq i\leq n+1$.
}
Under this condition, there exist two operators $\PR_{0}$ and $\RR_{0}$ acting on $\BR^{n+1}$ such that $\PR_{0}\RR_{0}=\RR_{0}\PR_{0}=\OR$, $\LR_{0}=\lam\PR_{0}+\RR_{0}$, $\PR_{0}$ is a projection onto the one-dimensional eigenspace of $\lam$, and the spectrum of $\RR_{0}$ does not intersect with small open ball in $\C$ at the center $\lam$.
Choose any $h_{0}\in \BR^{n+1}$ so that $\PR_{0}h_{0}\neq 0$ and let $h=\PR_{0}h_{0}/\|\PR_{0}h_{0}\|_{\BR^{n+1}}$. Then $\LR_{0}h=\lam h$. Since $\BR^{1}\supset \BR^{n+1}$ and $h$ is not zero, we see that $h$ is the corresponding eigenfunction of the eigenvalue $\lam$ of $\LR_{0}\,:\,\BR^{1}\to \BR^{1}$. Since $\lam$ is simple, for any $f\in \BR^{1}$, there exists a unique number $\nu(f)\in \C$ such that the projection $\PR_{0,1}\,:\,\BR^{1}\to \BR^{1}$ onto the eigenspace of $\lam$ has the form $\PR_{0,1}(f)=\nu(f)h$. Then $\nu$ is a bounded linear functional from $\BR^{1}$ to $\C$ with $\|\nu\|_{\BR^{1}\to \C}\leq \|\PR_{0,1}\|_{\BR^{1}}$ and $\nu(h)=1$. In summary, $(\lam,h,\nu)\in \C\times \BR^{n+1}\times (\BR^{1})^{*}$ satisfies
\alil
{
\LR_{0}h=\lam h, \ \LR_{0}^{*}\nu=\lam\nu,\ \nu(h)=1,\ \|h\|_{\BR^{n+1}}=1.\label{eq:L0_triple}
}
We give the eigenvalue $\lam(\e)$ and the corresponding eigenprojection $\PR(\e,\cd)\,:\,\BR^{n+1}\to \BR^{n+1}$ for $\LR(\e,\cd)$ such that $\lam(\e)\to \lam$ and $\|\PR(\e,\cd)-\PR_{0}\|_{\BR^{n+1}\to \BR^{0}}\to 0$ as $\e\to 0$. By a similar argument above, we get a triplet $(\lam(\e),h(\e,\cd),\nu(\e,\cd))\in \C\times \BR^{n+1}\times (\BR^{1})^{*}$ satisfying
\alil
{
\LR(\e,h(\e,\cd))=\lam(\e)h(\e,\cd), \ \LR(\e,\cd)^{*}\nu(\e,\cd)=\lam(\e)\nu(\e,\cd),\ \nu(\e,h(\e,\cd))=1,\ \|h(\e,\cd)\|_{\BR^{n+1}}=1\label{eq:Lt_triple}
}
for each $\e$.
In addition to the above condition (A.1), we introduce the following condition:
\ite
{
\item[(A.2)] The eigenvectors $\nu$ and $\nu(\e,\cd)$ given in (\ref{eq:L0_triple}) and (\ref{eq:Lt_triple}) are extended to bounded linear functionals on $\BR^{0}$, and satisfy that $\LR^{*}\nu=\lam\nu$, $\LR(\e,\cd)^{*}\nu(\e,\cd)=\lam(\e)\nu(\e,\cd)$ and $\limsup_{\e\to 0}\|\nu(\e,\cd)\|_{\BR^{0}\to \C}<\infty$.
}
Then we obtain the following:
\thm
{\label{th:GLtoGAPT}
Assume that (GL.1)-(GL.5) and (A.1)(A.2) are satisfied. Then the conditions (I)-(IV), the conditions in Theorem \ref{th:asymp_eval_nB_gene2} and the conditions in Theorem \ref{th:asymp_eval_nB_gene3} are all fulfilled.
}
\pros
The conditions (GL.4) and (GL.5) yield the condition (I). In particular, the property that $\LR_{0}$ is a linear operator acting on $\BR^{i}$ for all $0\leq i\leq n+1$ is guaranteed from (GL.5) letting $j=0$ and $\e=0$. 
The conditions (II) follows immediately from (A.1). To check the validity of (III) and (IV), we show $\nu(\e,h)\to 1$ and $\nu(h(\e,\cd))\to 1$ as $\e\to 0$. By $\|\PR(\e,\cd)-\PR_{0}\|_{\BR^{n+1}\to \BR^{0}}\to 0$, we have $\|\PR(\e,h)-h\|_{\BR^{0}}\to 0$. 
Furthermore, By
\ali
{
|\nu(h(\e,\cd))\nu(\e,h)-1|= |\nu(\PR(\e,h))-\nu(h)|\leq \|\nu\|_{\BR^{0}\to \K}\|\PR(\e,h)-h\|_{\BR^{0}}\to 0,
}
$|\nu(h(\e,\cd))|\leq \|\nu\|_{\BR^{0}\to \K}<\infty$ and $|\nu(\e,h)|\leq \sup_{\e>0}\|\PR(\e,\cd)\|_{\BR^{n+1}\to \BR^{0}}\|h\|_{\BR^{n+1}}<\infty$, we get $\nu(h(\e,\cd))\to 1$ and $\nu(\e,h)\to 1$. Therefore the conditions (III) and (IV) are valid. 
\smallskip
\par
Next we will check the conditions of Theorem \ref{th:asymp_eval_nB_gene2} and of Theorem \ref{th:asymp_eval_nB_gene3}.
The boundedness of $\|\kappa(\e,\cd)\|_{\BR^{0}\to \K}$ is yielded by the assumption $\|\nu(\e,\cd)\|_{\BR^{0}\to \K}$ in (A.2) and $\nu(\e,h)\to 1$. Moreover, remark the equation $\Delta_{j}(\e,\cd)=\tLR_{j-1}(\e,\cd)\e^{j-1}=\LR_{j}\e^{j}+\tLR_{j}(\e,\cd)\e^{j}$. Then (GL.5) says that 
\alil
{
&\|\tLR_{j-1}(\e,\cd)\|_{\BR^{i}\to \BR^{i-j}}\leq c_{\adr{btLRj}}\e\quad (1\leq j\leq n+1 \text{ and }j\leq i\leq n+1)\label{eq:Lj-1<=e}\\
&\|\tLR_{j}(\e,\cd)\|_{\BR^{i}\to \BR^{i-j}}\leq \|\LR_{j}\|_{\BR^{i}\to \BR^{i-j}}+c_{\adr{GLc5}}\quad (1\leq j\leq n+1 \text{ and }j\leq i\leq n+1)\label{eq:btLj}\\
&\|\LR(\e,\cd)\|_{\BR^{i}\to \BR^{i}}\leq c_{\adr{bLe}}\quad (0\leq i\leq n+1)\label{eq:bLe}
}
for some constants $c_{\adl{btLRj}}, c_{\adl{bLe}}>0$. In particular, $\|\tLR_{j}(\e,\cd)\|_{\BR^{j+1}\to \BR^{0}}\to 0$ for each $0\leq j\leq n$.  Thus the conditions of Theorem \ref{th:asymp_eval_nB_gene2} and Theorem \ref{th:asymp_eval_nB_gene3} are fulfilled.
\smallskip
\par
Finally we check the conditions of Theorem \ref{th:asymp_efunc_nB_gene}. The condition (a) of this theorem is clear. We will check the condition (b). The boundedness of $\nu$, $\LR_{j}$ and $\tLR_{j}(\e,\cd)$ immediately follow from (A.1), (GL.4) and (\ref{eq:btLj}), respectively. 
To prove the boundedness of $g(\e,\cd)$, we recall the form $g(\e,\cd)=h+\ti{g}(\e,\cd)$ with $\ti{g}_{0}(\e,\cd)=\ti{\lam}_{0}(\e)\SR g(\e,\cd)-\SR\tLR_{0}(\e,g(\e,\cd))$ in (\ref{eq:ge-h=...}). Observe that $\|g(\e,\cd)\|_{\BR^{n+1}}=1/\nu(g(\e,\cd))$ is bounded uniformly $\e>0$. Furthermore, we have that for $i=0,1,\dots, n$
\alil
{
|\|g(\e,\cd)\|_{\BR^{i}}-\|h\|_{\BR^{i}}|\leq&\|\ti{g}_{0}(\e,\cd)\|_{\BR^{i}}
\leq&
\|\SR\|_{\BR^{i}}(|\ti{\lam}_{0}(\e)|\|g(\e,\cd)\|_{\BR^{i}}+c_{\adr{btLRj}}\e\|g(\e,\cd)\|_{\BR^{i+1}}).\label{eq:ge-hinBi<=...}
}
If $i=n$ then the above inequality means $\sup_{\e>0}\|g(\e,\cd)\|_{\BR^{n}}<+\infty$. Indeed, if we suppose $\|g(\e,\cd)\|_{\BR^{n}}\to +\infty$ as $\e\to 0$, then the inequality yields $|1-\|h\|_{\BR^{n}}/\|g(\e,\cd)\|_{\BR^{n}}|\leq \|\SR\|_{\BR^{n}}o(1)$. This contradicts with $1\neq 0$. Thus $g(\e,\cd)$ is bounded in $\BR^{n}$. By using repeatedly (\ref{eq:ge-hinBi<=...}) for $i=n-1,n-2,\dots, 0$, we obtain $\sup_{\e>0}\|g(\e,\cd)\|_{\BR^{i}}<+\infty$ for any $i=0,1,\dots, n$.
\smallskip
\par
Since $\lam$ is in the resolvent of $\RR\,:\,\BR^{i}\to \BR^{i}$ for $i=1,2,\dots, n+1$, the resolvent $\RR_{\lam}$ is bounded on $\BR^{i}$. We see
\alil
{
&\|(\LR(\e,\cd)-\lam\IR)^{-1}-\RR_{\lam}\|_{\BR^{1}\to \BR^{0}}\leq c_{\adr{dGL}}\e^{\eta}\\
&\|(\LR(\e,\cd)-\lam\IR)^{-1}f\|_{\BR^{0}}\leq c_{\adr{bKL1}}\tau(\e)^{\eta}\|f\|_{\BR^{1}}+\tau(\e)^{-\eta}c_{\adr{bKL2}}\|f\|_{\BR^{0}}\label{eq:Le-lam^-1<=...}
}
for any small $\e>0$ and $f\in \BR^{1}$ for some constants $c_{\adl{dGL}},\eta, c_{\adl{bKL1}},c_{\adl{bKL2}}>0$, where the former follows from a result of \cite[Theorem 8.1]{GL} (or \cite{KL}) and the later is obtained by the estimate \cite[p.149, l.2]{KL} with $\tau(\e):=\|\tLR_{0}(\e,\cd)\|_{\BR^{1}\to \BR^{0}}$. These mean $\|\RR_{\lam}\|_{\BR^{1}\to \BR^{0}}<\infty$ in particular. To end the checking of the condition (b), we show that $\RR_{\lam}\,:\,\BR^{1}\to \BR^{0}$ is weak bounded.
Choose any $\e_{0}>0$ and take $\e>0$ so that $c_{\adr{bKL1}}\eta(\e)+c_{\adr{dGL}}\e^{\eta}\leq \e_{0}$. Then (\ref{eq:Le-lam^-1<=...}) implies that for any $f\in \BR^{1}$
\ali
{
\|\RR_{\lam}f\|_{\BR^{0}}\leq &\|(\LR(\e,\cd)-\lam\IR)^{-1}f-\RR_{\lam}f\|_{\BR^{0}}+\|(\LR(\e,\cd)-\lam\IR)^{-1}f\|_{\BR^{0}}\\
\leq&\tau(\e)^{-\eta}c_{\adr{bKL2}}\|f\|_{\BR^{0}}+\epsilon_{0}\|f\|_{\BR^{1}}.
}
This says that $\RR_{\lam}\,:\,\BR^{1}\to \BR^{0}$ is weak bounded. Thus the condition (b) is valid. Finally, the condition (c) is yielded by (\ref{eq:Lj-1<=e}).
 Consequently all conditions of Theorem \ref{th:asymp_eval_nB_gene3} are fulfilled. Hence the proof is complete.
\proe
\rem
{
The conditions (GL.2)(GL.3) are not necessary to obtain the asymptotic expansions of $\lam(\e),\nu(\e,\cd)$ under the conditions (A.1)(A.2). On the other hand, To get the expansion of the eigenfunction $h(\e,\cd)$, the conditions (GL.2)(GL.3) are used implicitly.
}
\rem
{
Let $X$ be $d$-dimensional $C^{\infty}$ compact connected Riemannian manifold, $T_{0}\in C^{r+1}(X,X)$ a transitive Anosov map with $r>1$, and $t\mapsto T_{t}\in C([0,1],C^{r+1}(X,X))$ perturbed Anosov maps. By virtue of this Theorem under suitable Banach spaces $\BR^{0},\cdots, \BR^{s}$ (see \cite[Section 9]{GL}), the Sinai-Ruelle-Bowen (SRB) measure $\mu(t,\cd)$ for $T_{t}$ has the asymptotic expansion with the coefficients and the remainder such as (\ref{eq:gk=2})(\ref{eq:tgke=...in2}) (see also \cite[Theorem 2.8]{GL}). In fact, the measure $\mu(t,\cd)$ is the corresponding eigenfunction of the simple eigenvalue $1$ of a transfer operator acting on $\BR^{0}$ for $T_{t}$ (see also \cite[Theorem 2.3]{GL}).
}
\appendix
\section{Thermodynamic formalism and Ruelle operators}\label{sec:thermo}
In this section, we recall the notions topological Markov shift and thermodynamic formalism treated in \cite{MU,Sar01,Sar03} which are used in our applications.
\smallskip
\par
Let $S$ be a countable set with distinct topology and $A=(A(ij))_{S\times S}$ a zero-one matrix. Consider the set
\alil
{
X=\{\om=\om_{0}\om_{1}\cdots\in\textstyle{\prod_{n=0}^{\infty}S}\,:\,A(\om_{k}\om_{k+1})=1 \text{ for any }k\geq 0\}\label{eq:X=...}
}
endowed with product topology induced by the discrete topology on $S$, and endowed with the shift transformation $\si\,:\,X\to X$ is defined by $(\si\om)_{n}=\om_{n+1}$ for any $n\geq 0$. This is called a {\it topological Markov shift} (or {\it countable Markov shift}) with state space $S$ and with transition zero-one matrix $A$. In what follows, we assume $X\neq \emptyset$.
An element $\om$ of $X$ denotes $\om=\om_{0}\om_{1}\om_{2}\cdots$ with $\om_{0},\om_{1},\om_{2},\dots \in S$.
Word $w=w_{1}w_{2}\dots w_{n}\in S^{n}$ is {\it admissible} if $A(w_{1}w_{2})=A(w_{2}w_{3})=\cdots=A(w_{n-1}w_{n})=1$.
For admissible word $w\in S^{n}$, {\it cylinder set} is defined by $[w]:=\{\om\in X\,:\,\om_{0}\cdots\om_{n-1}=w\}$. The transition matrix $A$ is {\it finitely irreducible} if there exists a finite subset $F\subset \bigcup_{n=1}^{\infty}S^{n}$ such that for any $a,b\in S$, there is $w\in F$ so that $a\cd w\cd b$ is admissible, where $a\cd w$ is the concatenation of $a$ and $w$. The matrix $A$ is called {\it finitely primitive} if there exists a finite subset $F\subset S^{N}$ with an integer $N\geq 1$ so that for any $a,b\in S$, $a\cd w\cd b$ is admissible for some $w\in F$ (see \cite{MU} for definition). Note that finitely primitively implies finitely irreducibility.
The matrix $A$ has the {\it big images and pre-images (BIP) property} if there is a finite set $S_{0}=\{a_{1},\cdots,a_{N}\}$ of $S$ such that for any $b\in S$, there exist $1\leq i,j\leq N$ such that $A(a_{i}b)=A(ba_{j})=1$.
It is known that the matrix $A$ is finitely irreducible if and only if $A$ is irreducible and fills the BIP property. Similarity, the matrix $A$ is finitely primitive if and only if $X$ is topologically mixing and $A$ has the BIP property.
\smallskip
\par
For $\theta\in (0,1)$, a metric $d_{\theta}\,:\,X\times X\to \R$ is defined by
$d_{\theta}(\om,\up)=\theta^{\min\{n\geq 0\,:\,\om_{n}\neq \up_{n}\}}$ if $\om\neq \up$ and $d_{\theta}(\om,\up)=0$ if $\om=\up$.
Then the metric topology induced by $d_{\theta}$ coincides with the product topology induced by the discrete topology on $S$.
Remark that $(X,d_{\theta})$ is a complete and separable metric space. Note also that if $\{a\in S\,:\,[a]\neq \emptyset\}$ is an infinite set, then $X$ is not compact.
\smallskip
\par
Next we introduce some function spaces and the notion of thermodynamic formalism \cite{MU,Sar06}. Put $\K=\R$ or $\C$. For integer $n\geq 1$ and function $f\,:\,X\to \K$, we let 
\ali
{
[f]_{\theta}^{n}:=\sup_{w\in S^{n}\,:\,[w]\neq \emptyset}\sup\{\frac{|f(\om)-f(\up)|}{d_{\theta}(\om,\up)}\,:\,\om,\up\in [w],\ \om\neq \up\}.
}
Note the inclusion $[f]_{\theta}^{n+1}\leq [f]_{\theta}^{n}$. 
A function $f\,:\,X\to \K$ is {\it $\theta$-locally Lipschitz continuous} if $[f]_{\theta}^{1}<\infty$. A function $f\,:\,X\to \K$ is called a {\it locally H\"older continuous} function if $f$ is a $\theta$-locally Lipschitz continuous for some $\theta\in (0,1)$.
We set
$C(X,\K)=\{f\,:\,X\to \K\,:\,\text{continuous functions}\}$, 
$C_{b}(X,\K)=\{f\in C(X,\K)\,:\,\|f\|_{C}<\infty\}$, 
$F_{\theta}(X,\K)=\{f\,:\,X\to \K\,:\,\theta\text{-locally Lipschitz continuous } \}$, and
$F_{\theta,b}(X,\K)=\{f\in F_{\theta}(X,\K)\,:\,\|f\|_{C}<\infty\}$ with
\ali
{
\|f\|_{F_{\theta}}:=\|f\|_{C}+[f]_{\theta}^{1}\quad\text{ and }\quad \|f\|_{C}:=\sup_{\om\in X}|f(\om)|.
}
Then the normed spaces $(C_{b}(X,\K),\|\cd\|_{C})$ and $(F_{\theta,b}(X,\K),\|\cd\|_{F_{\theta}})$ are Banach spaces. Remark the inclusion $F_{\theta}(X,\K)\subset F_{\theta^\p}(X,\K)$ for $\theta<\theta^\p$.
The symbol $\K$ may be omitted from these definitions when $\K=\C$.
For function $\ph\,:\,X\to \R$, the {\it topological pressure} $P(\ph)$ of $\ph$ is given by
\alil
{
P(\ph)=\lim_{n\to \infty}\frac{1}{n}\log \sum_{w\in S^{n}\,:\,[w]\neq \emptyset}\exp(\sup_{\om\in [w]}\sum_{k=0}^{n-1}\ph(\si^{k}\om))\label{eq:P(ph)=}
}
if it exists \cite{MU}.
If the transition matrix of $X$ is finitely irreducible and $\ph$ is in $F_{\theta}(X,\R)$ with $P(\ph)<\infty$, then $P(\ph)$ coincides with the Gurevich pressure $P_{G}(\ph)$ of $\ph$ which is introduced in \cite{Sar01}.
\smallskip
\par
A $\si$-invariant Borel probability measure $\mu$ on $X$ is said to be a {\it Gibbs measure} of the potential $\ph\,:\,X\to \R$ if there exist $c\geq 1$ and $P\in \R$ such that for any $\om\in X$ and $n\geq 1$
\ali
{
c^{-1}\leq \frac{\mu([\om_{0}\om_{1}\dots \om_{n-1}])}{\exp(-nP+\sum_{k=0}^{n-1}\ph(\si^{k}\om))}\leq c.
}
If $A$ is finitely irreducible and $\ph$ is in $F_{\theta}(X,\R)$  with $P(\ph)<\infty$, then the Gibbs measure of $\ph$ uniquely exists and $P$ equals $P(\ph)$ \cite{MU}.
\smallskip
\par
We end this section with the following result for transfer operators. 
For a real-valued function $\ph$ on $X$, the Ruelle operator $\LR_{\ph}$ associated to $\ph$ is defined by
\ali
{
\LR_{\ph} f(\om)=\sum_{a\in S\,:\,A(a\om_{0})=1}e^{\ph(a\cd \om)}f(a\cd\om)
}
if this series converges in $\C$ for a complex-valued function $f$ on $X$ and for $\om\in X$. It is known that if the incidence matrix is finitely irreducible and $\ph$ is in $F_{\theta}(X,\R)$ with finite topological pressure, then $\LR_{\ph}$ becomes a bounded linear operator both on $C_{b}(X)$ and $F_{\theta,b}(X)$. We state a version of Ruelle-Perron-Frobenius Theorem as follows.
\thm
{
\label{th:Spec_gap}
Let $X$ be a countable Markov shift with finitely irreducible transition matrix and $p$ the period of this matrix. Let $\ph\in F_{\theta}(X,\R)$ with finite pressure. Then the Ruelle operator $\LR_{\ph}\,:\,F_{\theta,b}(X)\to F_{\theta,b}(X)$ has the spectral decomposition
\ali
{
\LR_{\ph}=\sum_{i=0}^{p-1}\lam_{i}\PR_{i} +\RR
}
such that the following (1)-(4) are satisfied:
\ite
{
\item[(1)] $\lam_{i}=\lam\exp(i 2\pi\sqrt{-1}/p)$ is a simple eigenvalue of $\LR_{\ph}\,:\,F_{\theta,b}(X)\to F_{\theta,b}(X)$ and $\lam$ is the spectral radius of this operator;
\item[(2)] $\PR_{i}$ is the projection onto the one-dimensional eigenspace of $\lam_{i}$ and has the form
\ali
{
\PR_{i} f=\int_{X}f h_{i}\,d\nu_{i},
}
where $h_{i}$ is the corresponding  eigenfunction of the eigenvalue $\lam_{i}$ and $\nu_{i}$ is the corresponding positive eigenvector of $\lam_{i}$ of the dual operator $\LR_{\ph}^{*}\,:\,F_{\theta,b}(X)^{*}\to F_{\theta,b}(X)^{*}$ with $\nu_{i}(h_{i})=1$. In particular, $h:=h_{0}$ satisfies $c^{-1}\leq h\leq c$ for some constant $c\geq 1$ and $\var_{k}\log h_{0}\leq \sum_{i=k+1}^{\infty}\var_{i}\ph$ for any $k\geq 1$, where $\var_{k}f:=\sup\{|f(\om)-f(\up)|\,:\,\om_{i}=\up_{i}\text{ for any }i=0,1,\dots, k-1\}$. Moreover, the measure $\nu:=\nu_{0}$ is a Borel probability measure on $X$;
\item[(3)] $\PR_{i}\PR_{j}=\PR_{i}\RR=\RR\PR_{i}=\O$ for $i\neq j$ and the spectral radius of $\RR\,:\,F_{\theta,b}(X)\to F_{\theta,b}(X)$ is less than $\lam$;
\item[(4)] $P(\ph)$ equals $\log \lam$ and $h\nu$ becomes the Gibbs measure of the potential $\ph$.
}
}
The proof of this theorem is due to \cite{Aaronson_Denker} (see also \cite{Sar09}). Note that it is not hard to extend from the finitely primitive case to the finitely irreducible case.

\endthebibliography
\end{document}